\documentclass[10pt,twoside]{amsart} 
\usepackage{pslatex}
\usepackage{amssymb}
\theoremstyle{plain}

\newtheorem{theorem}{Theorem}[section]
\newtheorem{lemma}[theorem]{Lemma}
\newtheorem{proposition}[theorem]{Proposition}

\newtheorem{definition}[theorem]{Definition}
\newtheorem{remark}[theorem]{Remark}
\newtheorem{conjecture}[theorem]{Conjecture}

\dedicatory{Dedicated to Anatole Katok on the occasion of his 60th birthday.}
\title{Weakly mixing group actions: a brief survey and an example}
\author[V. Bergelson, A. Gorodnik]{V. Bergelson$^{*}$, A. Gorodnik}
\thanks{$^{*}$The first author was supported by the NSF Grants DMS-0070566 and DMS-0345350.}

\def\t{{}^t\!}
\def\plim{\mathop{p\hbox{-}\lim\,}}
\def\N{\mathbb N}

\def\Z{{\mathbb Z}}
\def\R{{\mathbb R}}

\begin{document}
\maketitle

\section{Introduction}

At its inception in the early 1930's, ergodic theory concerned itself with
continuous one-parameter flows of measure preserving transformations
(\cite{birk}, \cite{neum}, \cite{kn}, \cite{hop1}, \cite{hop2}). Soon it was realized that
working with $\mathbb{Z}$-actions rather than with $\mathbb{R}$-actions,
has certain advantages. On the one hand, while the proofs become simpler, the results
for $\mathbb{R}$-actions can often be easily derived from those for $\mathbb{Z}$-actions
(see, for example, \cite{kolm}). On the other hand, dealing with $\mathbb{Z}$-(or even with $\mathbb{N}$-) actions
extends the range of applications to measure preserving transformations
which are not necessarily embeddable in a flow.
Weakly mixing systems were introduced (under the name {\it
dynamical systems of continuous spectra}) in \cite{kn}.
By the time of publishing in 1937 of Hopf's book \cite{hop3},
the equivalence of the following conditions (which, for convenience, we formulate
for $\mathbb{Z}$-actions) was already known.
It is perhaps worth noticing that, while in most books either (i) or (ii)
below is taken as the ``official'' definition of weak mixing, the original definition in \cite{kn} corresponds to the condition (vi).

\begin{theorem}\label{th_mix_cond}
Let $T$ be an invertible measure-preserving
transformation of a probability measure space $(X, {\mathcal B}, \mu)$. Let
$U_T$ denote the operator defined on the space of measurable functions by $(U_Tf)(x) =
f(Tx)$. The following conditions are equivalent:
\begin{enumerate}
\item[(i)] For any $A,B\in {\mathcal B}$,
$$
\lim_{N\to\infty} \frac{1}{N} \sum_{n=0}^{N-1} \vert \mu(A\cap T^{-n}B) -
\mu(A)\mu(B)\vert = 0.
$$

\item[(ii)] For any $A,B\in {\mathcal B}$, there is a set $P\subset \N$ of density zero
such that
$$
\lim_{n\to\infty,\, n\notin P}   \mu(A\cap T^{-n}B) =  \mu(A)\mu(B).
$$

\item[(iii)] $T\times T$ is ergodic on the Cartesian square of $(X, {\mathcal B},
\mu)$.

\item[(iv)] For any ergodic probability measure preserving system $(Y,{\mathcal D},
\nu, S)$, the transformation $T\times S$ is ergodic on $X\times Y$.

\item[(v)] If $f$ is a measurable function such that for some $\lambda \in
{\mathbb C},\; U_Tf = \lambda f$ a.e., then $f= const$ a.e.

\item[(vi)] For   $f\in L^2(X, {\mathcal B}, \mu)$ with $\int_X f d\mu=0$, consider the
representation of the positive definite sequence $\langle U_T^n f, f\rangle , n\in \Z$,
as a Fourier transform of a measure $\nu$ on ${\mathbb T}=\mathbb{R}/\mathbb{Z}$:
$$
 \langle U_T^n f, f\rangle  = \int_ {\mathbb T} e^{2\pi inx} d\nu,\quad n\in \Z
$$
(this representation is guaranteed by Herglotz theorem, see \cite{he}). Then
$\nu$ has no atoms.
\end{enumerate}
\end{theorem}

\begin{remark}\label{12}
{\rm
It is not too hard to show that condition (i) can be replaced by the following
more general condition:
\begin{enumerate}
\item[($\hbox{i}^\prime$)] For any $A,B\in {\mathcal B}$ and any sequence of intervals $I_N=[a_N+1,a_N+2,\ldots,b_N]\subset\mathbb{Z}$, $N\ge 1$, with $|I_N|=b_N-a_N\to\infty$,
one has
$$
\lim_{N\to\infty} \frac{1}{|I_N|} \sum_{n=a_N+1}^{b_N} \vert \mu(A\cap T^{-n}B) -
\mu(A)\mu(B)\vert = 0.
$$
\end{enumerate}
Condition ($\hbox{i}^\prime$), in its turn, is equivalent to a still more general condition in which the
sequence of intervals $\{I_N\}_{N\ge 1}$ is replaced by an arbitrary {\it F\o lner}
sequence, i.e. a sequence of finite sets $F_N\subset \mathbb{Z}$, $N\ge 1$, such that
for any $a\in Z$,
$$
\frac{|(F_N+a)\cap F_N|}{|F_N|}\to 1\quad\hbox{as}\quad N\to\infty.
$$
This more general form of condition ($\hbox{i}^\prime$) makes sense for any (countably infinite)
amenable group and, as we shall see below (cf. Theorem \ref{th_dye}), can be used to
define the notion of weak mixing for actions of amenable groups.
} 
\end{remark}

\begin{remark}\label{r12}
{\rm 
If $(X, {\mathcal B}, \mu)$ is a {\it separable} space (which will be tacitly assumed from now on),
the condition (ii) can be replaced by the following condition (see Theorem I in \cite{kn}):
\begin{enumerate}
\item[($\hbox{ii}^\prime$)] There exists a set $P\subset\mathbb{N}$ of density zero such that for any $A,B\in {\mathcal B}$,
one has
$$
\mathop{\lim}_{n\to\infty,n\notin P} \mu (A\cap T^{-n} B)=\mu (A)\mu(B).
$$
\end{enumerate}
}
\end{remark}

Condition (ii) in Theorem \ref{th_mix_cond} indicates the subtle but
significant difference between weak and strong mixing: while
for strong mixing one has $\mu (A\cap T^{-n}B)\to \mu (A)\mu(B)$ as $n\to\pm\infty$
for {\it any} pair of measurable sets, a weakly mixing system which is not
strongly mixing is characterized by the {\it absence} of mixing for {\it some}
sets along {\it some} rarefied (i.e. having density zero) sequence of times.
Although the first examples of weakly but not strongly mixing
measure preserving transformations were quite complicated, numerous
classes of measure preserving systems that satisfy this property are
known by now.
For instance, one can show that the so-called interval exchange transformations (IET) are
often weakly mixing (\cite{ks}, \cite{ve}). On the other hand, A.~Katok
proved in \cite{ka} that the IET are never strongly mixing.
It should be also mentioned here that weakly mixing measure preserving transformations
are ``typical'', whereas strongly mixing ones are not
(see, for example, \cite{halm}). Before moving our discussion to weakly
mixing actions of general groups, we would like to formulate some more
recent results which exhibit new interesting facets of the notion of weak mixing.

\begin{theorem}\label{th_mix_cond2}
Let $T$ be an invertible measure-preserving
transformation of a probability measure space $(X, {\mathcal B}, \mu)$. The following conditions are equivalent:
\begin{enumerate}
\item[(i)] The transformation $T$ is weakly mixing.
\item[(ii)] Weakly independent sets are dense in ${\mathcal B}$. (Here a set $A\in {\mathcal
B}$ is {\rm weakly independent} if there exists a sequence $n_1 <
n_2 < \cdots $ such that the sets $T^{-n_i}A$, $i\ge 1$, are mutually
independent).
\item[(iii)] For any $A\in {\mathcal  B}$ and $k\in \N,\ k\ge 2,$ one has
$$
\lim_{N\to\infty} \frac{1}{N} \sum_{n=0}^{N-1} \mu(A\cap T^{-n}A\cap
T^{-2n}A\cap\cdots \cap T^{-kn}A ) = (\mu(A))^{k+1}.
$$
\item[(iv)] For any  $k\in \N,\ k\ge 2,$ any $f_1, f_2,...,f_k \in L^\infty(X,
{\mathcal B}, \mu)$, and any non-constant polynomials $p_1(n), p_2(n),...,
p_k(n) \in \Z [n]$ such that for all $i\neq j$, 
\hbox{$\mathrm{deg}(p_i - p_j) > 0$}, one has
$$
\lim_{N\to\infty} \frac{1}{N} \sum_{n=0}^{N-1} f_1(T^{p_1(n)}x)
f_2(T^{p_2(n)}x)\cdots  f_k(T^{p_k(n)}x) = \int\!\! f_1d\mu   \int\!\!
f_2d\mu\cdots \int\!\! f_kd\mu
$$
in $L^2$-norm.
\end{enumerate}
\end{theorem}

\begin{remark}
{\rm
Condition (ii) is due to U.~Krengel (see \cite{k} for this and related results).
Condition (iii) plays a crucial role in Furstenberg's ergodic proof of Szemer\'edi's
theorem on arithmetic progressions (see \cite{f_n} and \cite{fbook}).
Criterion (iv) was obtained in \cite{b_n}. Similarly to the ``linear'' case (iii),
the condition (iv) (or, actually, some variations of it) plays an important role
in proofs of polynomial extensions of Szemer\'edi's theorem
(see \cite{bl_n}, \cite{bm1}, \cite{bm2}, \cite{l}).
Note that the assumption $k\ge 2$ in (iii) and (iv) is essential.
Indeed, for $k=1$ condition (ii) expresses just the ergodicity of $T$, whereas for $k=1$, condition (iv)
is equivalent to the assertion that all non-zero powers of $T$ are ergodic.
The following equivalent form of condition (iv) is, however, both true and nontrivial already for $k=1$
(cf. condition ($\hbox{ii}^\prime$) in Remark \ref{r12}):
\begin{enumerate}
\item[($\hbox{iv}^\prime$)] 
For any $k\ge 1$ and any nonconstant polynomials $p_1(n),\ldots,p_k(n)\in\mathbb{Z}[n]$ such that for all $i\ne j$,
$\hbox{deg}(p_i-p_j)>0$, there exists a set $P\subset\mathbb{N}$ having zero density such that for any sets
$A_0,\ldots,A_k\in {\mathcal B}$, one has
$$
\mathop{\lim}_{n\to\infty,n\notin P} \mu (A_0\cap T^{p_1(n)} A_1\cap \cdots \cap T^{p_k(n)}A_k)=\mu (A_0)\mu(A_1)\ldots\mu(A_k).
$$
\end{enumerate}

}
\end{remark}

Theorems \ref{th_mix_cond}, \ref{th_mix_cond2}, and numerous appearances and applications of
weakly mixing one-parameter actions in ergodic theory hint that the notion of weak mixing could be of interest and of importance for actions of more general groups. One wants, of course, not only to be able to come up with a definition (this is not too hard: for example, condition (iii)
in Theorem \ref{th_mix_cond} makes sense for any group action), but also to be
able to have, similarly to the case of one-parameter actions, many diverse equivalent forms
of weak mixing including those which pertain to independence and higher degree
mixing properties of the type given in Theorem \ref{th_mix_cond2}.

Let $(T_g)_{g\in G}$ be a measure preserving action of a locally compact group $G$
on a probability measure space $(X, {\mathcal B}, \mu)$. If $G$ is amenable,
one can replace condition (i) in Theorem \ref{th_mix_cond} (or, rather, condition ($\hbox{i}^\prime$) in remark \ref{12})
by the assertion that the averages of the expressions
$|\mu(A\cap T_g B)-\mu(A)\mu(B)|$ taken along any F\o lner sequence in $G$
converge to zero. If $G$ is noncommutative, one also has to replace
condition (v) by the assertion that the
only finite-dimensional
subrepresentation of $(U_g)_{g\in G}$ (where $U_g$ is defined by
$(U_gf)(x)=f(T_g^{-1}x)$, $f\in L^2(X, {\mathcal B}, \mu)$)
is the restriction to the subspace of constant functions.
H.~Dye has shown in \cite{dye} that under these modifications
the conditions (i), (iii), and (v) in Theorem \ref{th_mix_cond}
are equivalent.
Dye's results are summarized in the following theorem
(cf. \cite{dye}, Corollary 1, p.~129). Again, for the sake of notational convenience,
we state the theorem for the case of a countable group $G$.

\begin{theorem}\label{th_dye}
Let $(T_g)_{g\in G}$ be a measure preserving action of a countable amenable group
$G$ on a probability measure space $(X, {\mathcal B}, \mu)$. Then
the following conditions are equivalent:
\begin{enumerate}
\item[(i)] For every F\o lner sequence $(F_n)_{n=1}^\infty$ in $G$ and any $A,B\in\mathcal{B}$,
one has
$$
\lim_{n\to\infty}\frac{1}{|F_n|}\sum_{g\in F_n} |\mu(A\cap T_g B)-\mu(A)\mu(B)|=0.
$$
\item[(ii)] The only finite dimensional subrepresentation of $(U_g)_{g\in G}$
is its restriction to the space of constant functions.
\item[(iii)] The diagonal action of $(T_g\times T_g)_{g\in G}$ on 
the product space $(X\times X, {\mathcal B}\otimes {\mathcal B}, \mu\otimes\mu)$
is ergodic (i.e. has no nontrivial invariant sets).
\end{enumerate}
\end{theorem}

\begin{remark} \label{r17}
{\rm
As a matter of fact, it is not too hard to show that conditions (ii) and (iii)
in Theorem \ref{th_dye} are equivalent for any locally compact noncompact
second countable group. See, for instance, \cite{mo}, Proposition~1, p.~157.
}
\end{remark}

A measure preserving system $(X, {\mathcal B}, \mu, T)$ is called a system with {\it discrete spectrum}
if $L^2(X, {\mathcal B}, \mu)$ is spanned by the eigenfunctions of the induced unitary operator $U_T$.
It is not hard to show that the condition (v) in Theorem  \ref{th_mix_cond} implies that a measure preserving
system $(X, {\mathcal B}, \mu, T)$ is weakly mixing if and only if it does not have a nontrivial 
factor which is a system with discrete spectrum. Remark \ref{r17} hints that a natural generalization
of this fact to general group actions holds as well.
(A measure preserving action of a group $G$ on a probability space $(X, {\mathcal B}, \mu)$  
has discrete spectrum if $L^2(X, {\mathcal B}, \mu)$ is representable as a direct sum of
finite-dimensional invariant subspaces.)

In \cite{neum2} and \cite{halm} von~Neumann and Halmos have shown that an ergodic one-parameter measure
preserving action has discrete spectrum if and only if it is conjugate to an action by rotations
on a compact abelian group. Again, this result has a natural extension to general group actions.
See \cite{mac} for details and further discussion.

The duality between the notion of weak mixing and discrete spectrum extends to the {\it relative case},
namely, to the situation where one studies the properties of a system relatively to its factors.
The theory of relative weak mixing is in the core of highly nontrivial structure theory
developed by H.~Furstenberg in the course of his proof (\cite{f_n}) of Szemer\'edi theorem.
See also \cite{fk} and \cite{fbook}, Chapter 6.

In \cite{zi1} and \cite{zi2} the duality between weak mixing and discrete spectrum is generalized to
extensions of general group actions. In particular, Zimmer established a far reaching
``relative'' version of Mackey's results on actions with discrete spectrum.




A useful interpretation of condition (i) in Theorem \ref{th_dye}
is that if $(T_g)_{g\in G}$ is a weakly mixing action of an amenable group
$G$, then for every $A,B\in\mathcal{B}$ and $\varepsilon>0$, the set
$$
R_{A,B}=\{g\in G: |\mu(A\cap T_g B)-\mu (A)\mu (B)|<\varepsilon\}
$$
is large in the sense that it has density $1$ with respect to any F\o lner
sequence $(F_n)_{n=1}^\infty$:
$$
\lim_{n\to\infty}\frac{|R\cap F_n|}{|F_n|}=1.
$$
A natural question that one is led to by this fact is whether there is a similar
characterization of the sets $R_{A,B}$ in the case when $G$ is not 
necessarily amenable.

It turns out that for every locally compact group which acts in a weakly mixing
fashion on a probability space, the set $R_{A,B}$ is always ``conull'',
and in more than one sense. One approach, undertaken in \cite{br},
is to utilize the classical fact that functions of the form $\psi(g)=\mu(A\cap T_g A)$
are positive definite. This implies that such $\psi (g)$, 
as well as a slightly more general functions of the form $\phi(g)=\mu(A\cap T_g B)$,
are {\it weakly almost periodic}
(see \cite{eb}).
By a theorem of Ryll-Nardzewski (see \cite{rn}), there is a unique invariant mean on
the space $\hbox{WAP}(G)$ of weakly almost periodic functions.
Denoting this mean by $M$ and assuming that for every $A,B\in\mathcal{B}$,
the function $g\mapsto \mu(A\cap T_gB)$ is continuous on $G$, let us call
the action $(T_g)_{g\in G}$ weakly mixing if for all 
$$
f_1,f_2\in L^2_0(X,\mathcal{B},\mu)\stackrel{def}{=}\{f\in L^2(X,\mathcal{B},\mu): \int_X fd\mu=0\},
$$
one has
$$
M\left(\left|\int_X f_1(x)f_2(T_gx)d\mu(x)\right|\right)=0.
$$

\begin{theorem}\label{th_br}
(\cite{br}, Theorem 4.1)
Let $(T_g)_{g\in G}$ be a measure preserving action of a locally compact second countable group $G$
on a probability space $(X,\mathcal{B},\mu)$.
The following are equivalent:
\begin{enumerate}
\item[(i)] $(T_g)_{g\in G}$ is weakly mixing.
\item[(ii)] For every $f_1,f_2\in L^2(X,\mathcal{B},\mu)$,
$$
M\left(\left|\int_X f_1(x)f_2(T_gx)d\mu(x)-\int f_1d\mu\int f_2d\mu\right|\right)=0.
$$
\item[(iii)] For every $f_0,\ldots,f_n\in L^2_0(X,\mathcal{B},\mu)$ and $\varepsilon >0$,
there exists $g\in G$ with
$$
\left|\int_X f_0(x)f_i(T_gx)d\mu(x)\right|<\varepsilon,\quad i=1,\ldots,n.
$$
\item[(iv)] For every $g_1,\ldots,g_n\in G$, $f\in L^2_0(X,\mathcal{B},\mu)$, and $\varepsilon >0$,
there exists $g\in G$ such that
$$
\left|\int_X f(T_gx)f(T_{g_i}x)d\mu(x)\right|<\varepsilon,\quad i=1,\ldots,n.
$$
\item[(v)] For all $F\in L^2(X,\mathcal{B},\mu)$, where $F$ is not equivalent to a constant, the set
$\{f(T_gx):g\in G\}$ is not relatively compact in $L^2(X,\mathcal{B},\mu)$.
\item[(vi)] $L^2_0(X,\mathcal{B},\mu)$ contains no nontrivial finite dimensional
invariant subspaces of $(U_g)_{g\in G}$.
\item[(vii)] $(T_g\times T_g)_{g\in G}$ is ergodic.
\item[(viii)] $(T_g\times T_g)_{g\in G}$ is weakly mixing.
\end{enumerate}
\end{theorem}

We shall describe now one more approach to weak mixing for general group actions
(see \cite{b}, Section~4, for more details and discussion).
Let $G$ be a countably infinite, not necessarily amenable discrete group.
For the purposes of the following discussion it will be convenient to view $\beta G$,
the Stone-\v{C}ech compactification of $G$, as the space of ultrafilters on $G$, i.e.
the space of $\{0,1\}$-valued finitely additive probability measures on the power set $\mathcal{P}(G)$ of $G$.
Since elements of $\beta G$ are $\{0,1\}$-valued measures, it is natural to identify
each $p\in\beta G$ with the set of all subsets having $p$-measure $1$,
and so we shall write $A\in p$ instead of $p(A)=1$.
(This explains the terminology: ultrafilters are just {\it maximal} filters.)
Given $p,q\in G$, one defines the product $p\cdot q$ by
$$
A\in p\cdot q \Leftrightarrow \{x:Ax^{-1}\in p\}\in q.
$$
The operation defined above is nothing but convolution of measures, which, on the other hand, is an extension of the group operation on $G$.
(Note that elements of $G$ are in one-to-one correspondence with point masses, the so-called
{\it principal} ultrafilters.)
It is not hard to check that the operation introduced above is associative and that $(\beta G,\cdot)$
is a left topological compact semigroup (which, alas, is never a group for infinite $G$).
For a comprehensive treatment of topological algebra in the Stone-\v{C}ech compactification,
the reader is referred to \cite{his}.
By a theorem due to R.~Ellis \cite{e1}, any compact semigroup with a left continuous
operation has an idempotent. (There are, actually, plenty of them since there are $2^\mathfrak{c}$
disjoint compact semigroups in $\beta G$.) Idempotent ultrafilters find numerous applications
in combinatorics (see, for example, \cite{hin_su} and \cite{his}, Part 3) and also are quite useful in ergodic
theory and topological dynamics (see, for example, \cite{b0}, \cite{b}).
Given an ultrafilter $p\in\beta G$ and a sequence $(x_g)_{g\in G}$ in a compact Hausdorff
space, one writes 
$$
\plim_{g\in G} x_g=y
$$
if for any neighborhood $U$ of $y$, one has
$$
\{g\in G: x_g\in U\}\in p.
$$
Note that in compact Hausdorff spaces p-limit always exists and is unique.

The following theorem, which is an ultrafilter analogue of Theorem 1.7 from \cite{fk2},
illustrates the natural connection between idempotents in $\beta G$
and ergodic theory of unitary actions.

\begin{theorem}\label{th_ult}
Let $(U_g)_{g\in G}$ be a unitary action of a countable group $G$ on a
Hilbert space~$\mathcal{H}$.
For any nonprincipal idempotent $p\in \beta G$ and any $f\in\mathcal{H}$ one has
$$
\plim_{g\in G} U_gf=Pf\quad (\hbox{weakly})
$$
where $P$ is the orthogonal projection on the subspace $\mathcal{H}_r$ of \hbox{\rm $p$-rigid elements},
that is, the space defined by
$$
\mathcal{H}_r=\{f:\plim U_gf=f\}.
$$
\end{theorem}

Theorem \ref{th_ult} has a strong resemblance to the classical von~Neumann's ergodic theorem.
In both theorems a generalized limit of $U_gf$, $g\in G$, (in case of von~Neumann's theorem
this is the Ces\'aro limit) is equal to an orthogonal projection of $f$ on a subspace of $\mathcal{H}$.
But while von~Neumann's theorem extends via Ces\'aro averages over F\o lner sets
to amenable groups only, Theorem \ref{th_ult} holds for nonamenable groups as well.

Given an element $p\in\beta G$, it is easy to see that $R=p\cdot\beta G$ is a right ideal in
$\beta G$ (that is, $R\cdot\beta G\subseteq R$). By using Zorn's lemma one can show
that any right ideal contains a minimal ideal. It is also not hard to prove that any
minimal right ideal in a compact left topological semigroup is closed
(see \cite{b}, Theorem 2.1 and Exercise 6). Now, by Ellis' theorem,
any minimal ideal in $\beta G$ contains an idempotent.
Idempotents belonging to minimal ideals are called minimal.
It is minimal idempotents which allow one to introduce a new characterization of weak mixing
for general groups. Recall that a set $A\subseteq \Z$ is called {\it syndetic}
if it has bounded gaps and {\it piecewise-syndetic} if it is an intersection of a syndetic set
with a union of arbitrarily long intervals. The following definition extends these notions
to general semigroups.

\begin{definition}
Let $G$ be a (discrete) semigroup.
\begin{enumerate}
\item[(i)] A set $A\subseteq G$ is called
syndetic if for some finite set $F\subset G$, one has 
$$
\bigcup_{t\in
F}At^{-1}=G.
$$
\item[(ii)] A set $A\subseteq G$ is piecewise syndetic if for some finite
set $F\subset G$, the family 
$$
\left\{ \left(\bigcup_{t\in F}At^{-1}\right)a^{-1}:\ a\in
G\right\}
$$
has the finite intersection property.
\end{enumerate}
\end{definition}

The following proposition establishes the connection between minimal idempotents and certain notions of largeness for subsets of $G$. It will be used below to give a new sense to the fact that for a weakly
mixing action on a probability space $(X,\mathcal{B},\mu)$, the set $R_{A,B}$ is large for all
$\varepsilon>0$ and $A,B\in\mathcal{B}$.

\begin{theorem}
(see \cite{b}, Exercise 7)
Let $G$ be a discrete semigroup and $p\in
(\beta G,\cdot)$ a minimal idempotent. Then
\begin{enumerate}
\item[(i)] For any $A\in p$, the set $B=\{g:\ Ag^{-1}\in p\}$ is syndetic.
\item[(ii)] Any $A\in p$ is piecewise syndetic.
\item[(iii)] For any $A\in p$, the set 
$$
A^{-1}A =
\{x\in G:yx\in A\hbox{ for some }y\in A\}
$$
is syndetic.  
(Note that if
$G$ is a group, then $A^{-1}A=
\{g_1^{-1}g_2 : g_1,g_2 \in A\}$.)
\end{enumerate}
\end{theorem}

\begin{definition}
A set $A\subseteq G$ is called central if there exists a minimal idempotent $p\in\beta G$ such that
$A\in p$. A set $A\subseteq G$ is called a $C^*$-set (or central$^*$ set) if $A$ is a member of any minimal
idempotent in $\beta G$.
\end{definition}

\begin{remark}
{\rm
The original definition of central sets (in $\Z$), which is due to Furstenberg (see \cite{fbook}, p.~161),
was the following: a subset $S\subseteq
{\mathbb N}$ is a central set if there exists a system $(X, T)$, a point $x\in X$,
a uniformly recurrent point $y$ proximal to $x$, and a neighborhood $U_y$
of $y$ such that $S=\{ n:\ T^nx\in U_y\}$.
The fact that central sets can be equivalently defined
as members of minimal idempotents was established in \cite{bh_nm}. See also Theorem 3.6 in \cite{b}.
}
\end{remark}

The following theorem gives yet another characterization of the
notion of weak mixing.

\begin{theorem}
(see \cite{b}, Section 4)
Let $(T_g)_{g\in G}$ be a measure preserving action of a countable group $G$ on a probability space
$(X,\mathcal{B},\mu)$. Then the following are equivalent:
\begin{enumerate}
\item[(i)] $(T_g)_{g\in G}$ is weakly mixing.
\item[(ii)] For every $f\in L^2(X,\mathcal{B},\mu)$ and any minimal idempotent $p\in\beta G$, one
has 
$$
\plim_{g\in G} f(T_gx)=\int_X fd\mu\quad\hbox{(weakly)}.
$$
\item[(iii)] There exists a minimal idempotent $p\in\beta G$ such that for any $f\in L^2(X,\mathcal{B},\mu)$,
one has $\plim_{g\in G} f(T_gx)=\int_X fd\mu$ (weakly).
\item[(iv)] For any $A,B\in {\mathcal B}$ and any $\varepsilon > 0$, the set
$$
\{ g\in G : \vert\mu(A\cap T_gB) - \mu(A)\mu(B)\vert  < \varepsilon\}
$$
is a C$^*$-set.
\end{enumerate}
\end{theorem}

Given a weakly mixing action of, say, a countable (but not necessarily amenable) group $G$, one would like
to know whether the action has higher order mixing properties along some massive and/or well-organized
subsets of $G$. For example, it is not hard to show that for any weakly mixing $\mathbb{Z}$-action
and any nonconstant polynomial $p(n)\in\mathbb{Z}[n]$, one can
find an $IP$-set $S$ such that for any $A,B \in\mathcal{B}$, one has 
$$
\mathop{\lim}_{n\to\infty,n\in S} \mu (A\cap T^{p(n)} B)=\mu (A)\mu(B).
$$
(An {\it $IP$}-set generated by a sequence $\{n_i:i\ge 1\}$
is, by definition, any set of the form
$\{n_{i_1}+\cdots+n_{i_k}:i_1<\cdots <i_k;\; k\in\mathbb{N}\}$.)

Another example of higher degree mixing along structured sets is provided by a theorem proved 
in \cite{brud}, according to which any weakly mixing action of a countable infinite direct sum $G=\oplus_{n\ge 1} \mathbb{Z}_p$,
where $\mathbb{Z}_p$ is the field of residues modulo $p$, has the property that the restriction of the action
of $G$ to an infinite subgroup (which is isomorphic to $G$) is Bernoulli
(see also \cite{bkm1}, \cite{bkm2}, \cite{bklm}, \cite{jrw}, \cite{j}, \cite{begun1}).

In Section 2 below we give a detailed analysis of higher order mixing properties for a concrete
classical example --- the standard action of $\hbox{\rm SL}(2,\mathbb{Z})$ on the $2$-dimensional torus $\mathbb{T}^2$.
Since $\hbox{\rm SL}(2,\mathbb{Z})$ contains mixing automorphisms (namely, hyperbolic automorphisms),
this action is weakly mixing. On the other hand, this action is not strongly mixing because 
$\hbox{\rm SL}(2,\mathbb{Z})$ contains nontrivial unipotent elements.

While many of the results obtained below hold (sometimes, after an appropriate modification)
for toral actions of $\hbox{\rm SL}(n,\mathbb{Z})$ and even in more general situations, we intentionally deal here
with $\hbox{\rm SL}(2,\mathbb{Z})$-actions in order to make the paper more accessible and important
issues more transparent.

Here is a sample of what is proved in the next section:
\begin{itemize}
\item 
(cf. Proposition \ref{eq_mix_main})
Let $T_1,\ldots, T_k\in\hbox{\rm SL}(2,\mathbb{Z})$. Then the following assertions are equivalent:
\begin{enumerate}
\item[(i)] For every $A_0,\ldots, A_k\in\mathcal{B}$,
$$
\mu(A_0\cap T_1^nA_1\cap \cdots \cap T_k^nA_k)\to \mu(A_0)\cdots\mu(A_k)\quad\hbox{as}\quad n\to\infty.
$$
\item[(ii)]
Each $T_i$ is hyperbolic, $T_i\ne \pm T_j$ for $i\ne j$, and for every $\rho>1$, there are at most
two matrices among $T_i$, $i=1,\ldots,k$, having an eigenvalue $\lambda$ such that $|\lambda|=\rho$.
\end{enumerate}

\item
(cf. Proposition \ref{p_rohlin})
Let $T_1,\ldots,T_k\in\hbox{\rm SL}(2,\mathbb{Z})$ be hyperbolic automorphisms. Denote by $\lambda_i$
the eigenvalue of $T_i$ such that $|\lambda_i|>1$. Put $a_{0,n}=0$, $n\ge 1$. Let $k\ge 1$ and $a_{i,n}\in\mathbb{Z}$, $i=1,\ldots,k$,
be such that
$$
\min\{ \left|\log |\lambda_i|\cdot a_{i,n}-\log |\lambda_j|\cdot a_{j,n}\right|:0\le i<j\le n\}\to\infty\quad\hbox{as}\quad n\to\infty.
$$
Then for every $A_0,\ldots, A_k\in\mathcal{B}$,
$$
\mu(A_0\cap T_1^{a_{1,n}}A_1\cap \cdots \cap T_k^{a_{k,n}}A_k)\to \mu(A_0)\cdots\mu(A_k)\quad\hbox{as}\quad n\to\infty.
$$
This result generalizes Rokhlin's theorem \cite{roh} in the case of $2$-dimensional torus.
See also Proposition \ref{p_um} for an analogue of this result for unipotent automorphisms.

\item While every abelian group of automorphisms $G$, which acts in a mixing fashion on $\mathbb{T}^2$, is mixing
of order $k$ for every $k\ge 1$, (that is,  for every $k\ge 1$ and sequences $g_{0,n}=e$, $g_{1,n},\ldots, g_{k,n}\in G$ such that
$$
g_{i,n}^{-1}g_{j,n}\to \infty\quad\hbox{as}\quad n\to\infty\quad\hbox{for}\quad 0\le i<j\le k,
$$
one has
$$
\mu(A_0\cap g_{1,n}A_1\cap \cdots \cap g_{k,n}A_k)\to \mu(A_0)\cdots\mu(A_k)\quad\hbox{as}\quad n\to\infty)
$$
a nonabelian group of automorphisms of $\mathbb{T}^2$ is never mixing of order $2$ (see Proposition \ref{p_no_3}).
Note that there are nonabelian groups of automorphisms that act in a mixing fashion on $\mathbb{T}^2$ (see 
the discussion after Proposition \ref{p-np}).
\end{itemize}

\section{$\hbox{\rm SL}(2,\mathbb{Z})$-action on torus}

\begin{definition}
A sequence $T_n\in\hbox{\rm SL}(2,\mathbb{Z})$, $n\ge 1$, is called mixing if for every
$f_1,f_2\in L^\infty (\mathbb{T}^2)$,
\begin{equation}\label{eq_mixing}
\int_{\mathbb{T}^2} f_1(T_n\xi)f_2(\xi)d\xi\to \left(\int_{\mathbb{T}^2}f_1(\xi)d\xi\right)\left(\int_{\mathbb{T}^2}f_2(\xi)d\xi\right)\quad \hbox{as}\quad n\to\infty.
\end{equation}
A transformation $T\in\hbox{\rm SL}(2,\mathbb{Z})$ is called mixing if
the sequence $T^n$, $n\ge 1$, is mixing.
\end{definition}
Note that this definition is different from the one given in \cite{berb}.

Recall that a matrix $T$ is called {\it hyperbolic} if its eigenvalues
have absolute values different from $1$, and {\it unipotent} if all its eigenvalues are equal to $1$.
It is well-known that an automorphism $T\in \hbox{\rm SL}(2,\mathbb{Z})$ is mixing on the torus
$\mathbb{T}^2$ if and only if it is hyperbolic. This implies that the action of $\hbox{\rm SL}(2,\mathbb{Z})$
on $\mathbb{T}^2$ is weakly but not strongly mixing and motivates the following problem:
give necessary and sufficient conditions for a sequence $T_n\in \hbox{\rm SL}(2,\mathbb{Z})$, $n\ge 1$,
to be mixing.

We start with a useful and straightforward lemma (cf. Theorem~3.1(1) in \cite{ber}).
For a matrix $T$, denote by $\t T$ its transpose.

\begin{lemma}\label{l_mix1}
A sequence $T_n\in \hbox{\rm SL}(2,\mathbb{Z})$, $n\ge 1$, is mixing if and only if for every
$(x,y)\in (\mathbb{Z}^2)^2-\{(0,0)\}$, the equality $\t T_n x+y=0$ holds for finitely many $n$ only.
\end{lemma}

\begin{proof}
To prove that $T_n$ is mixing, it is sufficient to check (\ref{eq_mixing}) for $f_1$ and $f_2$ in the dense
subspace of trigonometric polynomials. It follows that $T_n$ is mixing if and only if (\ref{eq_mixing}) holds for
$f_1$ and $f_2$ that are characters of the form
\begin{equation}\label{eq_xi}
\chi_{_x}(\xi)=e^{2\pi i\left<x,\xi\right>},\quad x\in\mathbb{Z}^2, \xi\in\mathbb{T}^2.
\end{equation}
For $x,y\in\mathbb{Z}^2$, one has
$$
\int_{\mathbb{T}^2} \chi_{_x}(T_n\xi)\chi_{_y}(\xi)d\xi=\int_{\mathbb{T}^2} \chi_{_{\t T_nx+y}}(\xi)d\xi=
\begin{cases}
0 & \hbox{if}\quad \t T_nx+y\ne 0,\\
1 & \hbox{if}\quad \t T_nx+y= 0.\\
\end{cases}
$$
It follows that for $(x,y)\in (\mathbb{Z}^2)^2-\{(0,0)\}$,
$$
\int_{\mathbb{T}^2} \chi_{_x}(T_n\xi)\chi_{_y}(\xi)d\xi\to \left(\int_{\mathbb{T}^2}\chi_{_x}(\xi)d\xi\right)\left(\int_{\mathbb{T}^2}\chi_{_y}(\xi)d\xi\right)=0\quad \hbox{as}\quad n\to\infty
$$
if and only if the equality $\t T_n x+y=0$ holds for finitely many $n$ only. This proves the lemma.
\end{proof}

Denote by $\hbox{\rm M}(2,\mathbb{K})$ the set of $2\times 2$-matrices over a field $\mathbb{K}$.
Using Lemma \ref{l_mix1}, we can now prove the following proposition.

\begin{proposition}\label{p_mix_1}
Let $T_n\in \hbox{\rm SL}(2,\mathbb{Z})$, $n\ge 1$, $\|\cdot\|$ be the $\max$-norm on $\hbox{\rm M}(2,\mathbb{R})$,
and $\mathcal{D}\subset\hbox{\rm M}(2,\mathbb{R})$ denote the set of limit points of the sequence $\frac{T_n}{\|T_n\|}$ as $n\to\infty$.
Then the sequence $T_n$ is not mixing if and only if there exist $A\in\hbox{\rm M}(2,\mathbb{Q})$ and $B\in\hbox{\rm M}(2,\mathbb{Q})$ such that
$B\in \mathcal{D}$ and $T_n=A+ \|T_n\|B$ for infinitely many $n\ge 1$.
\end{proposition}

\begin{proof}
We may assume that $\|T_n\|\rightarrow\infty$. (Indeed, 
if $\|T_n\|\nrightarrow\infty$, then there exists a matrix $T_0$ such that $T_n=T_0$ for infinitely many $n$, and the statement is obvious.)

``$\Leftarrow$'': Let $T_n=A+ \|T_n\|B$. Since
$$
\det B= \mathop{\lim}_{n\to\infty} \det \left(\frac{T_n}{\|T_n\|}\right)=\mathop{\lim}_{n\to\infty}\frac{1}{\|T_n\|^2}=0,
$$
$B$ is degenerate. Thus, there exists $x\in\mathbb{Z}^2-\{0\}$ such that $\;\t Bx=0$.
Then for infinitely many $n$, $\t T_nx=\t Ax$, and, by Lemma \ref{l_mix1}, $T_n$ is not mixing.

``$\Rightarrow$'': By Lemma \ref{l_mix1}, there exists $(x,y)\in (\mathbb{Z}^2)^2-\{(0,0)\}$
such that $\t T_nx=-y$ for infinitely many $n$.
By passing, if needed, to a subsequence,
we may assume that this equality holds for all $n\ge 1$.
It is clear that $\gcd (x_1,x_2)=\gcd (y_1,y_2)$.
Thus, we may assume that $x$ and $y$ are {\it primitive} (that is, the $\gcd$
of their coordinates is $1$). Take $C,D\in\hbox{\rm SL}(2,\mathbb{Z})$
such that $Ce_1=x$ and $De_1=-y$ where $e_1=(1,0)$. Then
$$
\t T_n=D\left(\begin{tabular}{rr}
$1$ & $a_n$\\
$0$ & $1$
\end{tabular}
\right)C^{-1}=DC^{-1}+a_n D\left(\begin{tabular}{rr}
$0$ & $1$\\
$0$ & $0$
\end{tabular}
\right)C^{-1}
$$
for some $a_n\in\mathbb{Z}$. Put $\t F_1=DC^{-1}$ and $\;\t F_2=D\left(\begin{tabular}{rr}
$0$ & $1$\\
$0$ & $0$
\end{tabular}
\right)C^{-1}$.
We have 
\begin{equation}\label{eq_m0}
T_n=F_1+a_nF_2,
\end{equation}
and
\begin{equation}\label{eq_m1}
|a_n|\cdot\|F_2\|-\|F_1\|\le \|T_n\|\le |a_n|\cdot\|F_2\|+\|F_1\|.
\end{equation}
Hence, $\|T_n\|\sim |a_n|\cdot\|F_2\|$ as $n\to\infty$.
Replacing, if necessary, $F_2$ by $-F_2$ and $a_n$ by $-a_n$ we may assume that $a_n>0$ for infinitely many $n$.
Then $B\stackrel{def}{=}\frac{F_2}{\|F_2\|}\in \mathcal{D}$. Passing to a subsequence, we get that $a_n>0$ for $n\ge 1$.
By triangle inequality and (\ref{eq_m1}),
$$
\left\| T_n- \|T_n\|\frac{F_2}{\|F_2\|} \right\| \le
\| T_n-a_nF_2 \| + \left\| a_n F_2-\|T_n\|\frac{F_2}{\|F_2\|} \right\|
=\|F_1\|+\left| a_n\|F_2\|-\|T_n\|\right|\le 2\|F_1\|.
$$
Thus, for infinitely many $n$, $T_n- \|T_n\|B=A$ for some $A\in \hbox{\rm M}(2,\mathbb{Q})$.
This proves the proposition.
\end{proof}

We illustrate the usefulness of Proposition \ref{p_mix_1} by the following two propositions.

\begin{proposition}
Let $U,V\in\hbox{\rm SL}(2,\mathbb{Z})$ be unipotent matrices.
Then the sequence $T_n=U^{-n}V^n$ is mixing if and only if $UV\ne VU$.
\end{proposition}

\begin{proof}
If $U$ and $V$ commute, one can show that they are powers of a single unipotent
transformation. Hence, in this case, the sequence $T_n=U^{-n}V^n$ is not mixing.

Conversely, suppose that $UV\ne VU$.
There exist $A,B\in\hbox{\rm SL}(2,\mathbb{Z})$ such that
$$
U=A^{-1}\left(\begin{tabular}{rr}
$1$ & $u$\\
$0$ & $1$
\end{tabular}
\right)A\quad\hbox{and}\quad
V=B^{-1}\left(\begin{tabular}{rr}
$1$ & $v$\\
$0$ & $1$
\end{tabular}
\right)B
$$
for some $u,v\in\mathbb{Z}-\{0\}$. It is sufficient to show that the sequence $S_n=A T_nB^{-1}$
is mixing. Let $AB^{-1}=\left(\begin{tabular}{rr}
$a$ & $b$\\
$c$ & $d$
\end{tabular}
\right)$. We have
$$
S_n=\left(\begin{tabular}{rr}
$1$ & $-nu$\\
$0$ & $1$
\end{tabular}
\right)AB^{-1}\left(\begin{tabular}{rr}
$1$ & $nv$\\
$0$ & $1$
\end{tabular}
\right)
=\left(\begin{tabular}{rr}
$a-(cu)n$ & $b-(av+du)n-(cv)n^2$\\
$c$ & $d+(cv)n$
\end{tabular}
\right).
$$
When $c=0$,
$$
V=B^{-1}\left(\begin{tabular}{rr}
$1$ & $v$\\
$0$ & $1$
\end{tabular}
\right)B
=B^{-1}(AB^{-1})^{-1}\left(\begin{tabular}{rr}
$1$ & $v$\\
$0$ & $1$
\end{tabular}
\right)(AB^{-1})B=
A^{-1}\left(\begin{tabular}{rr}
$1$ & $v$\\
$0$ & $1$
\end{tabular}
\right)A,
$$
and it follows that $U$ and $V$ commute. Thus, $c\ne 0$.

We apply now Proposition \ref{p_mix_1}. For sufficiently large $n$,
$\|S_n\|=|b-(av+du)n-(cv)n^2|$.
Also 
$$
\frac{S_n}{\|S_n\|}\longrightarrow\left(\begin{tabular}{rr}
$0$ & $-\hbox{sign}(cv)$\\
$0$ & $0$
\end{tabular}
\right)\stackrel{def}{=} C.
$$
Since $S_n-\|S_n\|C$ is not constant for infinitely many $n$,
the sequence $S_n$ is mixing.
\end{proof}

\begin{remark}
{\rm
When $U,V\in\hbox{\rm SL}(2,\mathbb{Z})$ are commuting unipotent transformations,
the sequence $U^{-n}V^n$ is relatively mixing in the sense of Definition \ref{def_rmix} below.

}
\end{remark}

Using a similar argument, one proves the following proposition:
\begin{proposition}
Let $U,V\in\hbox{\rm SL}(2,\mathbb{Z})$ such that $U$ is unipotent, and $V$ is hyperbolic.
Then the sequence $T_n=U^{-n}V^n$ is mixing.
\end{proposition}

\begin{proof}
Denote by $E_{ij}$ the $2\times 2$ matrix with $1$ in position $(i,j)$ and $0$'s elsewhere.
For some $A,B\in\hbox{GL}(2,\mathbb{R})$ and $\lambda$ with $|\lambda|>1$,
$$
U=A^{-1}\left(\begin{tabular}{cc}
$\lambda$ & $0$\\
$0$ & $\lambda^{-1}$
\end{tabular}
\right)A\quad\hbox{and}\quad
V=B^{-1}\left(\begin{tabular}{rr}
$1$ & $1$\\
$0$ & $1$
\end{tabular}
\right)B.
$$
We write
$$
T_n=A^{-1}\left(\begin{tabular}{cc}
$\lambda^{-n}$ & $0$\\
$0$ & $\lambda^n$
\end{tabular}
\right)AB^{-1}\left(\begin{tabular}{rr}
$1$ & $n$\\
$0$ & $1$
\end{tabular}
\right)B
=\lambda^{-n} C+ \lambda^{-n}n D+ \lambda^{n} E +\lambda^{n}n F
$$
where
$$
C=A^{-1}E_{11}A,\;\; D=A^{-1}E_{11}AB^{-1}E_{12}B,\;\; E=A^{-1}E_{22}A,\;\; F=A^{-1}E_{22}AB^{-1}E_{12}B.
$$

Suppose that $F\ne 0$. Then $\frac{T_n}{\|T_n\|}\to \frac{F}{\|F\|}$. By Proposition \ref{p_mix_1},
we need to show that there is no $X\in\hbox{M}(2,\mathbb{R})$ such that $T_n-\|T_n\|\frac{F}{\|F\|}=X$ for infinitely many $n$.
Since $F$ is degenerate, one of the matrices $C$, $D$, $E$ is not a scalar multiple of $F$ (say $C$).
Take a basis of $\hbox{\rm M}(2,\mathbb{R})$ which contains $C$ and $F$. The $C$-coordinate of
$T_n-\|T_n\|\frac{F}{\|F\|}$ with respect to this basis is equal to $\lambda^{-n}+\alpha\lambda^{-n}n+\beta \lambda^{n}$
for some $\alpha,\beta\in\mathbb{R}$. This shows that the sequence $T_n-\|T_n\|\frac{F}{\|F\|}$ 
consists of distinct matrices for sufficiently large $n$. Thus, $T_n$ is mixing.

Suppose that $F=0$. Then $\frac{T_n}{\|T_n\|}\to \frac{E}{\|E\|}$. By the same argument as in the
previous paragraph, the sequence $T_n-\|T_n\|\frac{E}{\|E\|}$ 
consists of distinct matrices for sufficiently large $n$. This implies that $T_n$ is mixing.
\end{proof}

\begin{remark}
{\rm
When $U,V\in\hbox{\rm SL}(2,\mathbb{Z})$ are hyperbolic and $U\ne V$, the sequence
$U^{-n}V^n$ is mixing. This follows from Proposition \ref{eq_mix_main} below.
}
\end{remark}

Next, we study multiple mixing for general sequences. 

\begin{definition}
Let $T_{i,n}\in \hbox{\rm SL}(2,\mathbb{Z})$, $n\ge 1$, i=1,\ldots,k.
The sequences $T_{1,n},\ldots,T_{k,n}$ are jointly mixing if for every $f_i\in L^\infty (\mathbb{T}^2)$, $i=1,\ldots,k+1$,
\begin{equation}\label{eq_mixm}
\int_{\mathbb{T}^2} f_1(T_{1,n}\xi)\cdots f_k(T_{k,n}\xi)f_{k+1}(\xi)d\xi\to \left(\int_{\mathbb{T}^2}f_1(\xi)d\xi\right)\cdots\left(\int_{\mathbb{T}^2}f_{k+1}(\xi)d\xi\right)
\end{equation}
as $n\to\infty$.
Transformations $T_1,\ldots,T_k$ are called jointly mixing if the sequences
$T_1^n,\ldots,T_k^n$, $n\ge 1$, are jointly mixing.
\end{definition}

In \cite{ber}, this property was called w-jointly strongly mixing (see Definition 3.6 in \cite{ber}).

In the course of proving Proposition \ref{eq_mix_main} below, we shall need the following
immediate extension of Lemma \ref{l_mix1} (cf. Theorem 4.3(1) in \cite{ber}).

\begin{lemma} \label{l_mixm}
Let $T_{i,n}\in \hbox{\rm SL}(2,\mathbb{Z})$, $n\ge 1$, i=1,\ldots,k.
The sequences $T_{1,n},\ldots,T_{k,n}$ are jointly mixing if and only if for every
$(x_1,\ldots, x_{k+1})\in (\mathbb{Z}^2)^{k+1}-\{(0,\ldots,0)\}$ the equality 
$$
\t T_{1,n} x_1+\cdots + \t T_{k,n}x_k+x_{k+1}=0
$$
holds for finitely many $n$ only.
\end{lemma}



\begin{proposition}\label{eq_mix_main}
Let $T_i\in\hbox{\rm SL}(2,\mathbb{Z})$, $i=1,\ldots,k$. The transformations $T_1,\ldots,T_k$
are jointly mixing if and only if each of $T_i$ is hyperbolic, $T_i\ne \pm T_j$ for $i\ne j$, and for every $\rho>1$, there are at most
two matrices among $T_i$, $i=1,\ldots,k$, having an eigenvalue $\lambda$ such that $|\lambda|=\rho$.
\end{proposition}

\begin{proof}
If a matrix $T\in\hbox{\rm SL}(2,\mathbb{Z})$ has complex eigenvalues, they are units in an imaginary
quadratic field. This implies that the eigenvalues of $T$ are roots of unity. Hence, the transformation $T$ is not mixing
on $\mathbb{T}^2$. Therefore, we may assume that the eigenvalues of $T$ are real.

Next, we note that one can assume without loss of generality that the eigenvalues of $T_i$, $i=1,\ldots,k$,
are positive. Indeed, put $\tilde{T}_i=-T_i$ if the eigenvalues of $T_i$ are negative and
$T_i$ otherwise. Clearly, transformations $T_i$, $i=1,\ldots, k$, are jointly
mixing if and only if transformations $\tilde{T}_i$, $i=1,\ldots, k$, are jointly mixing.

Let the transformations $T_1,\ldots,T_k$ be jointly mixing. Then each of
the sequences $T_i^n$ and $T_i^{-n} T_j^n$, $i\ne j$, is mixing too.
This implies that all $T_i$ are hyperbolic
and $T_i\ne T_j$ for $i\ne j$. To show that the conditions of the theorem are necessary, we consider 
transformations $T_1,T_2,T_3\in\hbox{\rm SL}(2,\mathbb{Z})$ that have the same eigenvalue $\lambda>1$.
We claim that there exists $(x_1,x_2,x_3)\in (\mathbb{Z}^2)^3-\{(0,0,0)\}$ such that
\begin{equation}\label{eq_non_mix}
\t T_1^n x_1+\t T_2^n x_2+\t T_3^n x_3=0
\end{equation}
for every $n\ge 1$,
which, in view of Lemma \ref{l_mixm}, implies that the sequences $T_1^n,T_2^n,T_3^n$ are not jointly mixing.

Since $T_i$, $i=1,2,3$, have the same eigenvalues, there exist $A,B\in\hbox{GL}(2,\mathbb{R})$
such that 
\begin{equation}\label{eq_T2}
T_2=A^{-1}T_1A\quad\hbox{and}\quad T_3=B^{-1}T_1B.
\end{equation}
Note that the matrix $A$ is a solution
of the matrix equation 
\begin{equation}\label{eq_lin}
XT_2=T_1X,
\end{equation}
which can be rewritten as a homogeneous system of linear equations
with rational coefficients. The set of rational solutions of (\ref{eq_lin}) is dense
in the space of real solutions. It follows that there exists a rational solution (\ref{eq_lin})
such that $\det(X)\ne 0$. This shows that we may choose $A$ and $B$ in $\hbox{GL}(2,\mathbb{Q})$.
For every $v\in\mathbb{R}^2$, $v=v_{+}+v_{-}$ where $v_{+}$ and $v_{-}$ are eigenvectors
of $\t T_1$ with eigenvalues $\lambda$ and $\lambda^{-1}$ respectively ($\lambda >1$).
Define linear maps $P_+:v\mapsto v_+$ and $P_-:v\mapsto v_-$. Then
\begin{equation}\label{eq_ppm}
\t T_1=\lambda P_+ +\lambda^{-1} P_-,\;\;\; P_+P_-=P_-P_+=0,\;\;\; P_\pm^2=P_\pm,\;\;\; P_++P_-=\hbox{id}.
\end{equation}
Note that $P_+,P_-\in\hbox{\rm M}(2,\mathbb{Q}(\sqrt{d}))$ for some $d\in\mathbb{N}$ determined by $\lambda$.
When $\sqrt{d}\in\mathbb{Q}$, $\lambda$ and $\lambda^{-1}$ are algebraic integers in $\mathbb{Q}$,
and it follows that that $\lambda=\pm 1$, which is a contradiction. Thus, $\sqrt{d}\notin\mathbb{Q}$.
Denote by $\sigma$ the nontrivial Galois automorphism of the field extension
$\mathbb{Q}(\sqrt{d})/\mathbb{Q}$. Then $\lambda^\sigma=\lambda^{-1}$ and
$(P_+)^{\sigma}=P_{-}$. Using (\ref{eq_T2}) and (\ref{eq_ppm}), we may rewrite equation (\ref{eq_non_mix}) as
$$
\lambda^n\left(P_+x_1+ \t A P_+ \t A^{-1}x_2+\t B P_+\t B^{-1}x_3\right)
+\lambda^{-n}\left(P_-x_1+ \t A P_- \t A^{-1}x_2+\t B P_-\t B^{-1}x_3\right)=0.
$$
The columns of the matrices $P_+$, $\t A P_+ \t A^{-1}$, and $\t B P_+\t B^{-1}$
lie in the vector space $\mathbb{Q}(\sqrt{d})^2$ that has dimension $4$ over $\mathbb{Q}$.
Thus, these columns are linearly dependent over $\mathbb{Q}$, and there exists
$(x_1,x_2,x_3)\in (\mathbb{Z}^2)^3-\{(0,0,0)\}$ such that
$$
P_+x_1+ \t A P_+ \t A^{-1}x_2+\t B P_+\t B^{-1}x_3=0.
$$
Applying $\sigma$ to this equality, we get
$$
P_-x_1+ \t A P_- \t A^{-1}x_2+\t B P_-\t B^{-1}x_3=0.
$$
This implies (\ref{eq_non_mix}) and proves that the conditions in the proposition are
necessary for mixing.

To prove sufficiency consider
$S_i,T_i\in\hbox{\rm SL}(2,\mathbb{Z})$, $i=1,\ldots,k$, such that $S_i$ and $T_i$
have the same eigenvalue $\lambda_i>1$, and $\lambda_i<\lambda_j$ for $i<j$.
We need to show that the transformations $S_1,T_1,\ldots,S_k,T_k$ are jointly mixing.
By Lemma \ref{l_mixm}, it is enough to prove that there is no
$(x_1,y_1,\ldots,x_k,y_k,z)\in(\mathbb{Z}^2)^{2k+1}-\{(0,\ldots,0)\}$
such that the equality
\begin{equation}\label{eq_no_mix_pf}
\t S_1^nx_1+\t T_1^ny_1+\cdots + \t S_k^nx_k+ \t T_k^ny_k+z=0.
\end{equation}
holds for infinitely many $n\ge 1$. Suppose that such a $(2k+1)$-tuple exists. Without loss of generality,
we may assume that $y_k\ne 0$. As above, we define $P_{i,\pm},Q_{i,\pm}\in\hbox{\rm M}(2,\mathbb{R})$
such that
\begin{eqnarray*}
\t S_i=\lambda_i P_{i,+} +\lambda_i^{-1} P_{i,-},\;\;\; P_{i,+}P_{i,-}=P_{i,-}P_{i,+}=0,\;\;\; P_{i,\pm}^2=P_{i,\pm},\;\;\; P_{i,+}+P_{i,-}=\hbox{id},\\
\t T_i=\lambda_i Q_{i,+} +\lambda_i^{-1} Q_{i,-},\;\;\; Q_{i,+}Q_{i,-}=Q_{i,-}Q_{i,+}=0,\;\;\; Q_{i,\pm}^2=Q_{i,\pm},\;\;\; Q_{i,+}+Q_{i,-}=\hbox{id}.
\end{eqnarray*}
Then (\ref{eq_no_mix_pf}) can be rewritten as
$$
\sum_{i=1}^k \lambda_{i}^{-n}(P_{i,-}x_i+Q_{i,-}y_i)+
\sum_{i=1}^k \lambda_{i}^{n}(P_{i,+}x_i+Q_{i,+}y_i)+z=0.
$$
Dividing this equality by $\lambda_k^n$ and taking a limit over a subsequence $n_j\to\infty$,
we deduce that 
$$
P_{k,+}x_k+Q_{k,+}y_k=0.
$$
Suppose that $Q_{k,+}y_k=0$. Then $y_k\ne 0$ is a rational eigenvector of $S_k$ with
eigenvalue $\lambda_k^{-1}$. It follows that $\lambda_k,\lambda_k^{-1}\in\mathbb{Q}$.
On the other hand, $\lambda_k$ and $\lambda_k^{-1}$ are algebraic integers. Hence, $\lambda_k=\pm 1$,
which is a contradiction. This shows that 
$$
v\stackrel{def}{=}P_{k,+}x_k=-Q_{k,+}y_k\ne 0.
$$
We have
$$
S_k v=\lambda_k v=T_k v.
$$
As above, we denote by $\sigma_k$ the nontrivial automorphism of the quadratic extension
$\mathbb{Q}(\lambda_k)/\mathbb{Q}$. Then
$$
P_{k,+}^{\sigma_k}=P_{k,-},\quad Q_{k,+}^{\sigma_k}=Q_{k,-},\quad \lambda_k^{\sigma_k}=\lambda_k^{-1},
$$
and it follows that
$$
S_k v^{\sigma_k}=\lambda_k^{-1} v^{\sigma_k}=T_k v^{\sigma_k}.
$$
Since $v$ and $v^{\sigma_k}$ are linearly independent, this implies that $S_k=T_k$,
which is a contradiction. Thus, (\ref{eq_no_mix_pf}) holds for finitely many $n$ only.
The proposition is proved.
\end{proof}

Proposition \ref{eq_mix_main} shows, in particular, that transformations $T_1,T_2,T_3\in\hbox{SL}(2,\Z)$
need not be jointly mixing even when every two of them are jointly mixing.
Nonetheless, pairwise conditions are sufficient to imply mixing in the commutative situation.

\begin{proposition}  \label{p_abel}
Let $T_1,\ldots,T_k\in\hbox{\rm SL}(2,\Z)$ be commuting automorphisms of 
$\mathbb{T}^2$. Then they are jointly mixing if and only if the transformations $T_i$ and $T_i^{-1}T_j$,
$i\ne j$, are mixing.
\end{proposition}

\begin{proof}
It is clear that if $T_1,\ldots,T_k$ are jointly mixing, then $T_i$ and $T_i^{-1}T_j$, $i\ne j$, are mixing.

Conversely, suppose that $T_i$ and $T_i^{-1}T_j$, $i\ne j$, are mixing. Then $T_i\ne \pm T_j$ for $i\ne j$.
Since $T_i$ and $T_j$ commute and are hyperbolic, they can be simultaneously reduced to the diagonal form.
Thus, if $T_i$ and $T_j$ have the same eigenvalues of the same modulus, then $T_i=\pm T_j^{\pm 1}$.
It follows that the conditions of Proposition \ref{eq_mix_main} are satisfied and hence, 
$T_1,\ldots,T_k$ are jointly mixing.
\end{proof}

One can show that the natural analog of Proposition \ref{p_abel} holds in every dimension.

In the case of a single measure preserving transformation $T$, $T$ is mixing if and only if
for every $k\ge 1$, the transformation $T^k$ is mixing.
In the following proposition, we investigate what happens for general sequences in our group:

\begin{proposition}\label{p212}
Let $T_n\in\hbox{\rm SL}(2,\mathbb{Z})$, $n\ge 1$, be hyperbolic automorphisms.
Let $\lambda_n$ be the eigenvalue of $T_n$ with $|\lambda_n|>1$.
\begin{enumerate}
\item For any $k\ge 1$, if the sequence $\lambda_n$ is bounded, then $T_n$ is mixing if and only if $T_n^k$ is mixing.
\item For any $k\ge 2$, if $\lambda_n\to\infty$, the sequence $T_n^k$ is always mixing.
\end{enumerate}
\end{proposition}

\begin{proof}
Let $t_n=\hbox{Trace}(T_n)$. Then $T_n$ is a root of its characteristic polynomial $x^2-t_nx+1$.
Using the polynomial identity:
$$
x^k=P(x)(x^2-t_nx+1)+\alpha_{n,k}x+\beta_{n,k}
$$
where $\alpha_{n,k},\beta_{n,k}\in \mathbb{Z}$,
$$
\alpha_{n,k}=\frac{\lambda_n^k-\lambda_n^{-k}}{\lambda_n-\lambda_n^{-1}},\quad
\beta_{n,k}=\frac{\lambda_n^{-k+1}-\lambda_n^{k-1}}{\lambda_n-\lambda_n^{-1}},
$$
we deduce that
\begin{equation}\label{eq_ab}
T_n^k=\alpha_{n,k}T_n+\beta_{n,k}.
\end{equation}

Suppose that $\lambda_n$, $n\ge 1$, is bounded. Then the sequences $\alpha_{n,k}$ and $\beta_{n,k}$ are bounded, hence 
take on only finitely many values. Hence, the equality
\begin{equation}\label{eq_tnk}
\t T_n^kx+y=\t T_n(\alpha_{n,k}x)+(\beta_{n,k}x+y)=0
\end{equation}
holds for some $(x,y)\in(\mathbb{Z}^2)^2-\{(0,0)\}$ and infinitely many $n$ if and only if the equality
$\t T_nx'+y'=0$ holds for some $(x',y')\in(\mathbb{Z}^2)^2-\{(0,0)\}$ and infinitely many $n$.
By Lemma \ref{l_mix1}, this proves the first part of the proposition.

We assume now that $\lambda_n\to\infty$. Then
$$
\alpha_{n,k}\sim\lambda_n^{k-1}\quad\hbox{and}\quad\beta_{n,k}\sim-\lambda_n^{k-2}\quad\hbox{as}\quad n\to\infty.
$$
By Lemma \ref{l_mix1}, it is sufficient to show that if (\ref{eq_tnk}) holds for infinitely many $n$,
then $x=y=0$. Suppose that (\ref{eq_tnk}) holds for infinitely many $n$. Dividing by $\alpha_{n,k}$
and taking a limit over a subsequence $n_j\to\infty$, we conclude that $T_{n_j}x\to 0$.
Since the sequence $T_{n_j}x$ is discrete, it follows that $x=0$, and $y=0$. Thus, $T_n^k$
is mixing.
\end{proof}

\begin{remark}
{\rm
Note that the statement in part (2) of Proposition \ref{p212} fails for $k=1$. For example, let
$$
T_n=\left(\begin{tabular}{cc}
$n$ & $n-1$\\
$1$ & $1$
\end{tabular}
\right),\quad n\ge 1,
$$
If $\lambda_n$ denotes the largest eigenvalue of $T_n$, then clearly, $\lambda_n\to\infty$.
However, the sequence $T_n$, $n\ge 1$, is not mixing.
(This follows from Proposition \ref{p_mix_1}.)
}
\end{remark}

Recall a theorem of Rokhlin \cite{roh}: 

\begin{theorem}
{\bf (Rokhlin)}
Let $T$ be a mixing automorphism of a compact abelian group.
Then the sequences $T^{a_{1,n}},\ldots,T^{a_{k,n}}$ are jointly mixing provided that
$$
\min\{|a_{i,n}-a_{j,n}|:0\le i<j\le n\}\to\infty\quad\hbox{as}\quad n\to\infty,
$$
where $a_{0,n}=0$.
\end{theorem}

The following proposition shows that a naive generalization
of Rokhlin's theorem to a general sequence of automorphisms $T_n$ is false. 

\begin{proposition}\label{p213}
Let $T_n\in\hbox{\rm SL}(2,\mathbb{Z})$, $n\ge 1$.
Denote by $\lambda_n$ the eigenvalue of $T_n$ such $|\lambda_n|>1$. If the sequence $\lambda_n$, $n\ge 1$,
is bounded, then for any choice of $a_i\in\mathbb{Z}$, $i=1,\ldots,k$, ($k>1$) the sequences
$T_n^{a_1},\ldots, T_n^{a_k}$ are not jointly mixing.
\end{proposition}

\begin{proof}
Without loss of generality, we may assume that $a_i>0$, $i=1,\ldots,k$.

By Lemma \ref{l_mixm}, it sufficient to show that there exists a tuple
$(x_1,\ldots,x_{k+1})\in(\mathbb{Z}^2)^{k+1}-\{(0,\ldots,0)\}$ such that for infinitely many $n$,
$$
\t T_n^{a_1}x_1+\cdots+\t T_n^{a_k}x_k+x_{k+1}=0.
$$
By (\ref{eq_tnk}), the last equality reduces to
\begin{equation}\label{eq_last}
\t T_n(\alpha_{n,a_1}x_1+\cdots+\alpha_{n,a_k}x_k)+\beta_{n,a_1}x_1+\cdots+\beta_{n,a_k}x_k+x_{k+1}=0.
\end{equation}
Since the sequence $\lambda_n$, $n\ge 1$, is bounded, 
the sequences $\alpha_{n,a_i}$ and $\beta_{n,a_i}$ are bounded too.
Thus, they are constant for infinitely many $n$. Now one can easily choose $x_i\in\mathbb{Z}^2$,
$i=1,\ldots,k+1$, not all zero, such that (\ref{eq_last}) holds.
For example, one can take all $x_i$'s to be multiples a fixed nonzero integer vector.
\end{proof}

\begin{remark} {\rm
Even if a sequence $T_n\in\hbox{\rm SL}(2,\mathbb{Z})$, $n\ge 1$, is such that
\begin{enumerate}
\item[(i)] $T_n$ is hyperbolic and mixing on $\mathbb{T}^2$ for all $n$,
\item[(ii)] $\lambda_n\to\infty$, where $\lambda_n$ is the eigenvalue of $T_n$ such that $\lambda_n>1$,
\end{enumerate}
the sequences $T_n$ and $T_n^2$ need not be jointly mixing.
For example, put
$$
T_n=\left(\begin{tabular}{cc}
$n$ & $n^2-1$\\
$1$ & $n$
\end{tabular}
\right),\quad n\ge 1.
$$
Then $\t T_n^2 x+\t T_ny+ z=0$ for $x=\t (0,1)$, $y=\t (-2,0)$, $z=\t (0,-2)$
which implies that the sequences $T_n$ and $T_n^2$ are not jointly mixing.
On the other hand, it follows from Proposition \ref{p_mix_1} that the sequence $T_n$ is mixing.
This example also demonstrates that pairwise conditions are not sufficient to guarantee
joint mixing even when the elements commute for every fixed $n$.
}
\end{remark}

We give here a generalization of Rokhlin's theorem in the case of the $2$-dimensional torus.
(A similar extension of Rokhlin's theorem holds in any dimension.)

\begin{proposition}\label{p_roh}
Let $T_n\in\hbox{\rm SL}(2,\mathbb{Z})$, $n\ge 1$.
Denote by $\lambda_n$ the eigenvalue of $T_n$ such that $|\lambda_n|\ge 1$.
Put $a_{0,n}=0$, $n\ge 1$. Let $k\ge 1$ and $a_{i,n}\in\mathbb{Z}$, $i=1,\ldots,k$.
Denote
$$
\gamma_n=\min\{|a_{i,n}-a_{j,n}|:0\le i<j\le n\}.
$$
Suppose that one of the following conditions holds:
\begin{enumerate}
\item The sequence $T_n$ is mixing, and 
$$
\left\{ \frac{\|T_n\|}{\lambda_n^{\gamma_n}}:\; n\ge 1\right\} \hbox{ is bounded}.
$$
\item $$\frac{\|T_n\|}{\lambda_n^{\gamma_n}}\to 0\quad\hbox{as}\quad n\to\infty.$$
\end{enumerate}
Then the sequences $T_n^{a_{1,n}},\ldots,T_n^{a_{k,n}}$ are jointly mixing.
\end{proposition}

\begin{remark}
{\rm
Part (2) of the theorem with $T_n=T$, $n\ge 1$, implies Rokhlin's theorem for the case of $2$-dimensional torus.
}
\end{remark}

\begin{proof}
Since $T_n$ is measure-preserving, we are allowed to replace $a_{i,n}$ by
$a_{i,n}-\min \{a_{i,n}:i=0,\ldots,k\}$. It follows that without loss of generality, we may assume that
$$
\min \{a_{i,n}:i=0,\ldots,k\}=0.
$$
Also by changing order and passing, if needed, to subsequences,
we may assume that 
$$
\max \{a_{i,n}:i=1,\ldots,k\}=a_{k,n}.
$$

Suppose that the sequences $T_n^{a_{1,n}},\ldots,T_n^{a_{k,n}}$ are not jointly mixing.
By Lemma \ref{l_mixm}, there exists a tuple $(x_1,\ldots,x_{k+1})\in(\mathbb{Z}^2)^{k+1}-\{(0,\ldots,0)\}$
such that the equality
$$
\t T_n^{a_{1,n}}x_1+\cdots+\t T_n^{a_{k,n}}x_k+x_{k+1}=0
$$
holds for infinitely many $n$. By (\ref{eq_ab}), the last equality is equivalent to
$$
\sum_{i=1}^k \alpha_{n,a_{i,n}}\t T_n x_i+ \sum_{i=1}^k \beta_{n,a_{i,n}} x_i   +x_{k+1}=0.
$$
Note that in both cases, $\lambda_n^{\gamma_n}\to\infty$ as $n\to\infty$. 
Therefore, it follows that $\lambda_n^{a_{n,k}-a_{n,i}}\to\infty$.
$$
\alpha_{n,a_{i,n}}\sim \frac{\lambda_n^{a_{i,n}}}{\lambda_n-\lambda_n^{-1}}, \quad
\beta_{n,a_{i,n}}\sim \frac{-\lambda_n^{a_{i,n}-1}}{\lambda_n-\lambda_n^{-1}},\quad i=1,\ldots,k.
$$
Then
\begin{eqnarray*}
\t T_n x_k &=& -\sum_{i=1}^{k-1} \frac{\alpha_{n,a_{i,n}}}{\alpha_{n,a_{k,n}}}\t T_n x_i- \sum_{i=1}^k \frac{\beta_{n,a_{i,n}}}{\alpha_{n,a_{k,n}}} x_i -\frac{x_{k+1}}{\alpha_{n,a_{k,n}}}\\
&=& O\left(\frac{\|T_n\|}{\lambda_n^{\gamma_n}}\right)+O(\lambda_n^{-\gamma_n})+\lambda_n^{-1}x_k.
\end{eqnarray*}

Assume that condition (1) holds. Then the sequence $\t T_nx_k$ is bounded by infinitely many $n$.
Thus, it is constant for infinitely many $n$. It follows from Lemma \ref{l_mix1} that $x_k=0$.

Suppose that condition (2) holds. We prove that $x_k=0$.
If $\lambda_n\to\infty$, then $\t T_n x_k\to 0$ as $n\to\infty$, and this implies that $x_k=0$.
Otherwise, the sequences $\lambda_n$ and $\t T_n x_k$ are bounded for infinitely many $n$,
and consequently, they are constant for infinitely many $n$. Thus, $\t T_{n_j} x_k=\lambda_{n_j}^{-1}x_k$
for a subsequence $n_j$, and if $x_k\ne 0$, then $\lambda_{n_j},\lambda_{n_j}^{-1}\in\mathbb{Q}$.
Since $\lambda_{n_j}$ is an algebraic integer, $\lambda_{n_j}=\pm 1$, which contradicts condition (2).
This shows that $x_k=0$.

Now the proof can be completed by induction on $k$.
\end{proof}

\begin{remark}
{\rm
Condition (1) is not necessary for joint mixing of the sequences $T_n^{a_{1,n}},\ldots,T_n^{a_{k,n}}$.
For example, put
$$
T_n=\left(\begin{tabular}{cc}
$n^2$ & $n^3-1$\\
$1$ & $n$
\end{tabular}
\right),\; n\ge 1,\quad\hbox{and}\quad a_{i,n}=i,\; i=1,2.
$$
Even though $\frac{\|T_n\|}{\lambda_n}\to\infty$ as $n\to\infty$,
one can check with the help of Lemma \ref{l_mixm} that the sequences $T_n$ and $T_n^2$ are jointly mixing.
It would be of interest to find a necessary and sufficient condition for joint mixing
of sequences of the form $T_n^{a_{1,n}},\ldots,T_n^{a_{k,n}}$.
}
\end{remark}

The following proposition is yet another generalization of Rokhlin's theorem.

\begin{proposition} \label{p_rohlin}
Let $T_1,\ldots,T_k\in\hbox{\rm SL}(2,\mathbb{Z})$ be hyperbolic automorphisms. Denote by $\lambda_i$
the eigenvalue of $T_i$ such that $|\lambda_i|>1$. Put $a_{0,n}=0$, $n\ge 1$. Let $k\ge 1$ and $a_{i,n}\in\mathbb{Z}$, $i=1,\ldots,k$
be such that
\begin{equation}\label{eq_roh}
\min\{ \left|\log |\lambda_i|\cdot a_{i,n}-\log |\lambda_j|\cdot a_{j,n}\right|:0\le i<j\le n\}\to\infty\quad\hbox{as}\quad n\to\infty.
\end{equation}
Then the sequences $T_1^{a_{1,n}},\ldots,T_k^{a_{k,n}}$ are jointly mixing.
\end{proposition}

\begin{proof}
As in the proof of Proposition \ref{p_roh}, we reduce the proof to the case when
$$
\log |\lambda_{i+1}|\cdot a_{i+1,n}-\log |\lambda_i|\cdot a_{i,n}\to\infty\quad\hbox{as}\quad n\to\infty
$$
for $i=0,\ldots,k-1$.

Suppose that the sequences $T_1^{a_{1,n}},\ldots,T_k^{a_{k,n}}$ are not jointly mixing. By Lemma \ref{l_mixm},
there exists $(x_1,\ldots,x_{k},y)\in(\mathbb{Z}^2)^{k+1}-\{(0,\ldots,0)\}$ such that the equality
\begin{equation}\label{eq_last0}
\t T_1^{a_{1,n}}x_1+\cdots+\t T_k^{a_{k,n}}x_k +y=0
\end{equation}
holds for infinitely many $n$. Let $P_{i,+},P_{i,-}\in\hbox{M}(2,\mathbb{R})$, $i=1,\ldots,k$, be such that
$$
\t T_i=\lambda_i P_{i,+} +\lambda_i^{-1} P_{i,-},\;\;\; P_{i,+}P_{i,-}=P_{i,-}P_{i,+}=0,\;\;\; P_{i,\pm}^2=P_{i,\pm},\;\;\; P_{i,+}+P_{i,-}=\hbox{id}.
$$
By (\ref{eq_last0}),
\begin{equation}\label{eq_last1}
\sum_{i=1}^k \lambda_i^{a_{i,n}} P_{i,+}x_i+\sum_{i=1}^k \lambda_i^{-a_{i,n}} P_{i,-}x_i+y=0
\end{equation}
holds by infinitely many $n$. Dividing by $\lambda_k^{a_{k,n}}$ and taking limit over a subsequence $n_j\to\infty$,
we conclude that $P_{k,+}x_k=0$. 

If $x_k\ne 0$, it is an eigenvector of $\t T_k$ with the eigenvalue $\lambda_k^{-1}$. This implies that
$\lambda_k\in\mathbb{Q}$. On the other hand, $\lambda_k$ is an algebraic integer. Thus, $\lambda_k=\pm 1$.
This contradiction shows that $x_k=0$.
Using induction on $k$, we deduce from (\ref{eq_last1}) that $x_i=0$ for $i=1,\ldots,k$.
This shows that the sequences $T_1^{a_{1,n}},\ldots,T_k^{a_{k,n}}$ are jointly mixing.
\end{proof}

\begin{remark}
{\rm
It clear that condition (\ref{eq_roh}) in Proposition \ref{p_rohlin} follows from the following condition:
$$
a_{1,n}\to\infty\quad\hbox{and}\quad \frac{a_{i+1,n}}{a_{i,n}}\to\infty,\;\; i=1,\ldots,k,\quad\hbox{as}\quad n\to\infty,
$$
which also appears in Proposition \ref{p_um}.
}
\end{remark}

\begin{definition}\label{def_rmix}
Let $T_{i}\in\hbox{\rm SL}(2,\mathbb{Z})$, $i=1,\ldots,k$. Denote by $P_i:L^2(\mathbb{T}^2)\to L^2(\mathbb{T}^2)$,
$i=1,\ldots,k$, the orthogonal projection on the subspace of $T_i$-invariant functions.
Let $a_{i,n}\in \mathbb{Z}$, $i=1,\ldots,k$, $n\ge 1$.
We call the sequences $T_1^{a_{1,n}},\ldots,T_k^{a_{k,n}}$ relatively jointly mixing if for every $f_i\in L^\infty (\mathbb{T}^2)$,
$i=1,\ldots,k+1$,
\begin{equation}\label{eq_rel_mixing}
\int_{\mathbb{T}^2} f_1(T_1^{a_{1,n}}\xi)\cdots f_k(T_k^{a_{k,n}}\xi)f_{k+1}(\xi)d\xi\to \int_{\mathbb{T}^2}(P_1f_1)(\xi)\cdots (P_kf_k)(\xi)f_{k+1}(\xi)d\xi
\quad \hbox{as}\quad n\to\infty.
\end{equation}
\end{definition}

We have the following criterion for relative joint mixing of tuples of unipotent elements:

\begin{proposition}\label{p_relm}
Let $T_i\in\hbox{\rm SL}(2,\mathbb{Z})$, $i=1,\ldots,k$, be unipotent elements. Denote by
$v_i$, $i=1,\ldots,k$, a nonzero vector such that $\t T_iv_i=v_i$. 
Let $a_{i,n}\in \mathbb{Z}$, $i=1,\ldots,k$, $n\ge 1$. 
Then the sequences $T_1^{a_{1,n}},\ldots,T_k^{a_{k,n}}$ are relatively jointly mixing
if and only if for every $(\alpha_1,\ldots,\alpha_k)\in\mathbb{Z}^k-\{(0,\ldots,0)\}$
and $z\in\mathbb{Z}^2-\{0\}$, the equality
\begin{equation}\label{eq_rem}
\sum_{i=1}^k \alpha_i a_{i,n} v_i+z=0
\end{equation}
holds for finitely many $n$ only.
\end{proposition}

\begin{proof}
For some $A_i\in\hbox{\rm SL}(2,\mathbb{Z})$ and $s_i\in\mathbb{Z}-\{0\}$,
\begin{equation}\label{eq_unip}
\t T_i^{a_{i,n}}=A_i^{-1}\left(\begin{tabular}{cc}
$1$ & $s_i a_{i,n}$\\
$0$ & $1$
\end{tabular}
\right)A_i=E+s_ia_{i,n}B_i
\end{equation}
where $B_i=A_i^{-1}\left(\begin{tabular}{rr}
$0$ & $1$\\
$0$ & $0$
\end{tabular}
\right)A_i\in\hbox{\rm SL}(2,\mathbb{Z})$, and $E$ is the identity matrix.
To establish relative mixing, it is sufficient to check (\ref{eq_rel_mixing}) in the case when
$f_i$, $i=1,\ldots,k$, are characters of the form (\ref{eq_xi}).
For $x_1,\ldots,x_{k+1}\in\mathbb{Z}^2$, one has
\begin{eqnarray*}
\int_{\mathbb{T}^2} \chi_{_{x_1}}(T_1^{a_{1,n}}\xi)\cdots \chi_{_{x_k}}(T_k^{a_{k,n}}\xi)\chi_{_{x_{k+1}}}(\xi)d\xi
&=&\int_{\mathbb{T}^2} \chi_{_{(\t T_1^{a_{1,n}}x_1+\cdots +\t T_k^{a_{k,n}}x_k+x_{k+1})}}(\xi)d\xi\\
&=&
\begin{cases}
1 & \hbox{if}\quad \t T_1^{a_{1,n}}x_1+\cdots+ \t T_k^{a_{k,n}}x_k+x_{k+1}=0,\\
0 & \hbox{if}\quad \t T_1^{a_{1,n}}x_1+\cdots+ \t T_k^{a_{k,n}}x_k+x_{k+1}\ne 0.\\
\end{cases}
\end{eqnarray*}
Note that for every $x\in\mathbb{Z}^2$,
$$
P_i\chi_{_x}=
\begin{cases}
\chi_{_x} & \hbox{if}\quad B_ix=0,\\
0 & \hbox{if}\quad B_ix\ne 0.\\
\end{cases}
$$
Thus, (\ref{eq_rel_mixing}) always holds for $f_i=\chi_{_{x_i}}$ provided that $B_ix_i=0$
for all $i=1,\ldots,k$. It follows that the sequences $T_1^{a_{1,n}},\ldots,T_k^{a_{k,n}}$
are relatively jointly mixing if and only if for every $(x_1,\ldots,x_k)
\in(\mathbb{Z}^2)^k$ such that for some $i=1,\ldots,k$,
$B_ix_i\ne 0$ (equivalently, $T_ix_i\ne x_i$) the equality
\begin{equation}\label{eq_eq}
\t T_1^{a_{1,n}}x_1+\cdots +\t T_k^{a_{k,n}}x_k+x_{k+1}=0
\end{equation}
holds for finitely many $n$ only. By (\ref{eq_unip}), the last equality is equivalent to
$$
\sum_{i=1}^k s_ia_{i,n}B_ix_i +z=0
$$
where $z=\sum_{i=1}^{k+1} x_i$. Note that the columns of the matrix $B_i$ are rational
multiples of the vector $v_i$. Thus, $s_iB_ix_i=\alpha_iv_i$ for some $\alpha_i\in\mathbb{Q}$.
Since $B_ix_i\ne 0$ for some $i$, $(\alpha_1,\ldots,\alpha_k)\ne (0,\ldots,0)$.
This shows that (\ref{eq_eq}) holds if and only if
$$
\sum_{i=1}^k \alpha_i a_{i,n} v_i+z=0
$$
for some  $(\alpha_1,\ldots,\alpha_k)\in \mathbb{Q}^k-\{(0,\ldots,0)\}$ and $z\in\mathbb{Z}$.
Multiplying by a fixed integer, we get that $\alpha_i\in\mathbb{Z}$.
This proves the proposition.
\end{proof}

We record here a convenient corollary of Proposition \ref{p_relm}.

\begin{proposition}\label{p_um}
Let $T_i\in\hbox{\rm SL}(2,\mathbb{Z})$, $i=1,\ldots,k$, be unipotent elements, and
$a_{i,n}\in \mathbb{Z}$, $i=1,\ldots,k$, $n\ge 1$, such that
$$
a_{1,n}\to\infty\quad\hbox{and}\quad \frac{a_{i+1,n}}{a_{i,n}}\to\infty,\;\; i=1,\ldots,k,\quad\hbox{as}\quad n\to\infty.
$$
Then the sequences $T_1^{a_{1,n}},\ldots,T_k^{a_{k,n}}$ are relatively jointly mixing.
\end{proposition}

\begin{proof}
Suppose that the sequences $T_1^{a_{1,n}},\ldots,T_k^{a_{k,n}}$ are not relatively jointly mixing.
Then by Proposition \ref{p_relm}, (\ref{eq_rem}) holds for infinitely many $n$.
Dividing (\ref{eq_rem}) by $a_{k,n}$ and taking the limit over a subsequence $n_s\to\infty$,
we deduce that $\alpha_k=0$.
Similarly, it follows that $\alpha_i=0$ for $i=1,\ldots,k$. This shows that
the sequences $T_1^{a_{1,n}},\ldots,T_k^{a_{k,n}}$ are relatively jointly mixing.
\end{proof}

Let $T,S\in\hbox{\rm SL}(d,\mathbb{Z})$.
It was observed by Boshernitzan that it follows from the fact that the set of common periodic points
of $T$ and $S$ is dense in $\mathbb{T}^d$ that
for every nonempty open subset $\mathcal{U}$ of $\mathbb{T}^d$,
$$
\mathcal{U}\cap T^n\mathcal{U}\cap S^n\mathcal{U}\ne \emptyset
$$
for infinitely many $n$.
A measurable analogue of this fact is far less trivial. The following conjecture
seems plausible:

\begin{conjecture}\label{c_rec}
Let $T,S\in\hbox{\rm SL}(d,\mathbb{Z})$, and let $\mathcal{D}$ be a Borel subset of $\mathbb{T}^d$
of positive measure. Then
$$
\limsup_{n\to\infty}\; \mu(\mathcal{D}\cap T^n\mathcal{D}\cap S^n\mathcal{D})>0.
$$
\end{conjecture}

In fact, in all known to us examples,
$$
\limsup_{n\to\infty}\; \mu(\mathcal{D}\cap T^n\mathcal{D}\cap S^n\mathcal{D})\ge \mu(\mathcal{D})^3.
$$

\begin{remark}
{\rm
Note that when $T$ and $S$ generate a (virtually) nilpotent group, Conjecture
\ref{c_rec} follows from a general ``nilpotent'' multiple recurrence theorem proved
in \cite{l} (see also Theorem E in \cite{bl4}). It was, however, proved in \cite{bl},
that for any finitely generated solvable group of exponential growth $G$,
there exist a measure preserving action $(T_g)_{g\in G}$ on a probability
space $(X,\mathcal{B},\mu)$, elements $g_1, g_2\in G$ and a set $D\in\mathcal{B}$ with $\mu(D) >0$
such that for $T=T_{g_1}$ and $S=T_{g_2}$, one has $\mu(D\cap T^n D \cap S^n D)=0$ for all $n\ne 0$.
Nevertheless, we believe that for our special action of $\hbox{SL}(d,\mathbb{Z})$ on $\mathbb{T}^d$,
the Conjecture is true.
}
\end{remark}

We obtain below some partial results on the conjecture in the case of the $2$-dimensional torus.
Note that when $T$ and $S$ are hyperbolic the conjecture follows from Proposition
\ref{eq_mix_main}. In fact, in this case,
$$
\lim_{n\to\infty}\; \mu(\mathcal{D}\cap T^n\mathcal{D}\cap S^n\mathcal{D})=
\begin{cases}
\mu(\mathcal{D})^2 & \hbox{if}\quad T=S,\\
\mu(\mathcal{D})^3 & \hbox{if}\quad T\ne \pm S\\
\end{cases}
$$
and when $T=-S$, the limit set of the sequence $\mu(\mathcal{D}\cap T^n\mathcal{D}\cap S^n\mathcal{D})$
consists of two numbers: $\mu(\mathcal{D})^2$, $\mu(\mathcal{D}\cap -\mathcal{D})\mu(\mathcal{D})$. In particular, this shows
that $\liminf$ might be $0$ even when $\mu(\mathcal{D})>0$.

We can also settle the case when $T$ and $S$ are unipotent and hyperbolic respectively.
For this, we need a lemma:

\begin{lemma}\label{l_rec_help}
Let $T\in\hbox{\rm SL}(2,\mathbb{Z})$ be unipotent, and $S\in\hbox{\rm SL}(2,\mathbb{Z})$ hyperbolic.
Then the sequences $T^n$ and $S^n$, $n\ge 1$, are relatively jointly mixing.
\end{lemma}

\begin{proof}
As in the proof of Proposition \ref{p_relm}, it is sufficient to show that for
every $x,y,z\in\mathbb{Z}^2$ such that either $Tx\ne x$ or $y\ne 0$, the equality
\begin{equation}\label{eq_con_mix}
\t T^nx+\t S^n y+z=0
\end{equation}
holds for finitely many $n$ only. We have 
$$
\t T^n=E+nB
$$
where $E$ is the identity matrix and $B\in\hbox{\rm M}(2,\mathbb{Z})$. Let $\lambda$ be the eigenvalue of $S$ such that $|\lambda|>1$.
For some $P_+,P_-\in\hbox{\rm M}(2,\mathbb{R})$,
$$
\t S^n=\lambda^n P_+ +\lambda^{-n} P_-,\;\;\; P_+P_-=P_-P_+=0,\;\;\; P_\pm^2=P_\pm,\;\;\; P_++P_-=\hbox{id}.
$$
Equality (\ref{eq_con_mix}) is equivalent to
$$
\lambda^n P_+y +\lambda^{-n} P_-y +nBx+(x+z)=0.
$$
Suppose that it holds for infinitely many $n$. Dividing by $\lambda^n$
and taking the limit as $n\to\infty$, we deduce that $P_+y=0$.
Then $y$ is an eigenvector of $S$. If $y\ne 0$, then $y$ is 
a rational eigenvector of $S$, and $\lambda$ and $\lambda^{-1}$
are rational numbers that are algebraic integers. Hence, $\lambda=\pm 1$,
which is a contradiction. This implies that $y=0$.
Then it follows that $Bx=0$ (equivalently, $Tx=x$).
This shows that (\ref{eq_con_mix}) holds for finitely many $n$ only.
Thus, the sequences $T^n$ and $S^n$, $n\ge 1$, are relatively jointly mixing.
\end{proof}

Lemma \ref{l_rec_help} implies the following special case of Conjecture \ref{c_rec}.

\begin{proposition}
Let $T\in\hbox{\rm SL}(2,\mathbb{Z})$ be unipotent, and $S\in\hbox{\rm SL}(2,\mathbb{Z})$ hyperbolic.
Then for any measurable $\mathcal{D}\subseteq \mathbb{T}^2$, the limit of 
$\mu(\mathcal{D}\cap T^n\mathcal{D}\cap S^n\mathcal{D})$ as $n\to\infty$ exists, and
$$
\lim_{n\to\infty}\; \mu(\mathcal{D}\cap T^n\mathcal{D}\cap S^n\mathcal{D})\ge\mu(\mathcal{D})^3.
$$
Moreover, the equality holds if and only if $\mu(\mathcal{D})=1\;\hbox{or}\;0$.
\end{proposition}

\begin{proof}
Let $f$ be the characteristic function of the set $\mathcal{D}$. Denote by $P_T$ and $P_S$
the orthogonal projections on the the spaces of $T$- and $S$-invariant functions
respectively. Since $S$ is ergodic, $P_S f=\mu(\mathcal{D})$. By Lemma \ref{l_rec_help},
$$
\lim_{n\to\infty}\; \mu(\mathcal{D}\cap T^n\mathcal{D}\cap S^n\mathcal{D})
=\int_{\mathbb{T}^2} f(P_Tf)(P_Sf)d\mu=\mu(\mathcal{D})\|P_Tf\|_2^2\ge\mu(\mathcal{D})^3.
$$
\end{proof}

In the case when $T$ and $S$ are unipotent, Conjecture \ref{c_rec} seems to be open in general.
We prove a partial result for sets of special form.
For a function $f\in L^\infty(\mathbb{T}^d)$, its Fourier coefficients are denoted by
$$
\hat f(x)=\int_{\mathbb{T}^d}f(\xi)\chi_{_{-x}}(\xi)d\xi,\quad x\in\mathbb{Z}^d.
$$

\begin{proposition}
Let $T,S\in\hbox{\rm SL}(2,\mathbb{Z})$ be unipotent.
\begin{enumerate}
\item For any measurable $\mathcal{D}\subseteq\mathbb{T}^2$, the limit 
of $\mu(\mathcal{D}\cap T^n\mathcal{D}\cap S^n\mathcal{D})$ as $n\to\infty$ exists.

\item Suppose that $TS\ne ST$.
Let $A\in\hbox{\rm SL}(2,\mathbb{Z})$ be such that $A^{-1}TA$ is lower triangular unipotent.
Then for every set of the form $\mathcal{D}=A(\mathcal{D}_1\times\mathcal{D}_2)$ where $\mathcal{D}_1$
and $\mathcal{D}_2$ are measurable subsets of $\mathbb{T}^1$, 
$$
\lim_{n\to\infty} \mu(\mathcal{D}\cap T^n\mathcal{D}\cap S^n\mathcal{D})\ge\mu(\mathcal{D})^3.
$$
Moreover, the equality holds if and only if $\mu(\mathcal{D})=1\;\hbox{or}\;0$.
\end{enumerate}
\end{proposition}

\begin{proof}
We prove (1) in the case when $T$ and $S$ do not commute. (When $T$ and $S$ commute, 
they are powers of the same transformation, and the proof goes along the same lines as the proof below.)

Let $v$ and $w$ be primitive integer vectors such that $\t Tv=v$ and $\t S w=w$. We claim that for $f,g,h\in C^\infty(\mathbb{T}^2)$,
\begin{equation}\label{eq_recur}
\int_{\mathbb{T}^2} f(\xi)g(T^n\xi)h(S^n\xi)d\xi\to \sum_{i,j\in\mathbb{Z}} \hat f(-iv-jw)\hat g(iv)\hat h(jw)\quad\hbox{as}\quad n\to\infty.
\end{equation}
It follows from a standard argument that it is sufficient to check (\ref{eq_recur}) when $f$, $g$, $h$
are characters of the form (\ref{eq_xi}).

Let $f=\chi_{_x}$, $g=\chi_{_y}$, and $h=\chi_{_z}$ for some $x,y,z\in\mathbb{Z}^2$.
First, suppose that $x=-iv-jw$, $y=iv$, $z=jw$ for some $i,j\in\mathbb{Z}^2$. Then
$$
\int_{\mathbb{T}^2} f(\xi)g(T^n\xi)h(S^n\xi)d\xi =\int_{\mathbb{T}^2} \chi_{_{x+\t T^ny+\t S^nz}}(\xi)d\xi=1.
$$
This implies (\ref{eq_recur}) in this case.

Now we consider the case when $x$, $y$, $z$ are not of the above form.
We need to show that the equality $x+\t T^ny+\t S^nz=0$ holds for finitely many $n$ only.
Suppose that it holds for infinitely many $n$. Write $\t T =E+B$ and $\t S=E+C$ where $E$ is the identity matrix and $B,C\in\hbox{\rm SL}(2,\mathbb{Z})$
such that $B^2=C^2=0$. Then
$$
x+\t T^ny+\t S^nz=(x+y+z)+n(By+Cz)=0.
$$
holds for infinitely many $n$. This implies that $x+y+z=0$ and $By=-Cz$.
Note that the columns of matrix $B$ are multiples of the vector $v$, and the columns of
$C$ are multiples of $w$. If $By\ne 0$, $v$ is multiple of $w$, and it follows that in 
some basis of $\mathbb{R}^2$ both $T$ and $S$ are unipotent upper triangular. Then $TS=ST$,
and this contradicts the initial assumption. Thus, $By=Cz=0$. Equivalently,
$Ty=y$ and $Sz=z$. Hence, $x=-iv-jv$, $y=iv$, and $z=jw$ for some $i,j\in\mathbb{Z}$. This is a contradiction.
We have proved (\ref{eq_recur}).

Replacing $T$ by $A^{-1}TA$ and $S$ by $A^{-1}SA$, we reduce the problem to the case when $T$ is
lower triangular and $\mathcal{D}=\mathcal{D}_1\times\mathcal{D}_2$. Then $v=\t (1,0)$. Let $w=\t (a,b)$.
Let $f$ be the characteristic function of the set $\mathcal{D}$, and
$f_1$ and $f_2$ be characteristic functions of the sets $\mathcal{D}_1$ and $\mathcal{D}_2$ respectively.
Note that for $s,t\in\mathbb{Z}$, $\hat f(s,t)=\hat f_1(s)\hat f_2(t)$.
To prove part (2), we need to show
that
$$
\sum_{i,j\in\mathbb{Z}} \hat f(-i-aj,-bj)\hat f(i,0)\hat f(aj,bj)\ge\mu(\mathcal{D})^3.
$$
Using the Plancherel formula and the fact that $f_1^2=f_1$, we have
\begin{eqnarray*}
\sum_{i,j\in\mathbb{Z}} \hat f(-i-aj,-bj)\hat f(i,0)\hat f(aj,bj)&=&
\sum_{j\in\mathbb{Z}} \left(\sum_{i\in\mathbb{Z}}\hat f_1(-aj-i)\hat f_1(i)\right)
\hat f_2(-bj)\hat f_2(0)\hat f(aj,bj)\\
&=&\sum_{j\in\mathbb{Z}} \hat{(f^2_1)}(-aj) \hat f_2(-bj)\mu(\mathcal{D}_2)\hat f(aj,bj)\\
&=&\mu(\mathcal{D}_2)\sum_{j\in\mathbb{Z}} |\hat f_1(aj)|^2 |\hat f_2(bj)|^2\\
&\ge&\mu(\mathcal{D}_2)|\hat f_1(0)|^2 |\hat f_2(0)|^2\ge \mu(\mathcal{D})^3.
\end{eqnarray*}
We are done.
\end{proof}

Next, we investigate mixing properties of subgroups of $\hbox{\rm SL}(2,\mathbb{Z})$.

\begin{proposition}\label{p-np}
Let $H$ be a subgroup of $\hbox{\rm SL}(2,\mathbb{Z})$.
The action of $H$ on $\mathbb{T}^2$ is mixing if and only if $H$ contains
no nontrivial unipotent elements.
\end{proposition}

\begin{proof}
If the action of $H$ is mixing, then the action of every infinite subgroup of $H$ is mixing,
and consequently, $H$ does not contain nontrivial unipotent elements.

Conversely, suppose that the action of $H$ is not mixing.
By Lemma \ref{l_mix1}, there exists a sequence $h_n\in H$, $n\ge 1$, and $(x,y)\in(\mathbb{Z}^2)^2-\{(0,0)\}$
such that $\t h_nx=-y$ for every $n\ge 1$ and $h_n\to\infty$. Then $\t h_1^{-1}\t h_n x=x$ 
for every $n\ge 1$. Thus, $h_nh_1^{-1}\in H$ is a nontrivial unipotent element for
sufficiently large $n$. This proves the proposition.
\end{proof}

A subgroup of $\hbox{\rm SL}(2,\mathbb{Z})$ is called {\it nonparabolic} if it contains no nontrivial
unipotent elements. Nonparabolic subgroups are of interest from the point of view of ergodic
theory because they are precisely the groups that act in a mixing fashion on the torus $\mathbb{T}^2$.
It follows from the pigeonhole principle that every subgroup of $\hbox{\rm SL}(2,\mathbb{Z})$ of finite index
contains a nontrivial unipotent element. First examples of nonparabolic subgroups were
constructed by B.~H.~Neumann in \cite{neu} (see also \cite{mag}). Any Neumann subgroup has the
property that powers of a single unipotent element form a complete
system of representatives of the cosets of this group. In particular, Neumann subgroups are
maximal nonparabolic subgroups. There are examples of maximal nonparabolic subgroups
that are not Neumann (see \cite{t}, \cite{bl1}, \cite{bl2}). If $F$ is a free normal subgroup of finite
index in $\hbox{\rm SL}(2,\mathbb{Z})$ which is not equal to the commutant of $\hbox{\rm SL}(2,\mathbb{Z})$,
then the commutant of $F$ is nonparabolic (see \cite{mas}).

Although there are large subgroups in $\hbox{\rm SL}(2,\mathbb{Z})$ (e.g.
Neumann subgroup) whose actions on the torus
$\mathbb{T}^2$ are mixing, the following proposition shows that for nonabelian subgroups the higher
order of mixing is impossible.

\begin{proposition}\label{p_no_3}
A nonabelian subgroup of $\hbox{\rm SL}(2,\mathbb{Z})$ cannot be mixing of order $2$.
\end{proposition}

\begin{proof}
Let $H$ be a nonabelian subgroup of $\hbox{\rm SL}(2,\mathbb{Z})$.
Suppose that $H$ is mixing of order $2$. Take $g,h\in H$ such that $gh\ne hg$.
Since $H$ is mixing, $g$ and $h$ are hyperbolic.
Note that if $g^2h= hg^2$, $g^2$ and $h$ can be both reduced to the diagonal form, and this
would imply that $g$ and $h$ commute. Thus, $g^2h\ne hg^2$.
Put $h_i=g^{-i}hg^{i}$, $i=1,2,3$. Comparing eigenvalues, we deduce that $h_i\ne -h_j$.
Also, $h_i\ne h_j$ for $i\ne j$, since otherwise it would follow that $g$ and $h$ commute.
Therefore, by Proposition \ref{eq_mix_main}, the transformations $h_i$ and $h_j$ are jointly mixing
for $i\ne j$. In particular, $h_i^{-n}h_j^n\to\infty$ and $h_i^n\to\infty$ as $n\to\infty$.
On the other hand, since $h_1$, $h_2$, $h_3$
have the same eigenvalues, it follows from Proposition \ref{eq_mix_main} (and its proof) that
for some $f_1,f_2,f_3\in L^\infty(\mathbb{T}^2)$,
$$
\lim_{n\to\infty}\int_{\mathbb{T}^2} f_1(h_1^n\xi)f_2(h_2^n\xi)f_3(h_3^n\xi)d\xi\neq \int_{\mathbb{T}^2} f_1(\xi)d\xi
\int_{\mathbb{T}^2} f_2(\xi)d\xi\int_{\mathbb{T}^2} f_3(\xi)d\xi.$$
This implies that the sequences $h_1^{-n}h_2^n$ and $h_1^{-n}h_3^n$ are not jointly mixing,
which, in its turn, contradicts the assumption that the group $H$ is mixing of order $2$.
\end{proof}

More generally, a similar argument allows one to show that the standard action on the
torus $\mathbb{T}^d$ of a subgroup $H$ of $\hbox{\rm SL}(d,\Z)$ which is not virtually
abelian can not be mixing of order $d$.
Another approach to the proof of this fact can be found in \cite{bha}
where it is utilized for derivation of isomorphism rigidity for the action of $H$.


\begin{remark}
{\rm
It should be noted that in contrast with the nonabelean situation, any mixing $\mathbb{Z}^d$-action
on $\mathbb{T}^d$ is mixing of all orders (see \cite[Corollary~27.7]{sch}).
On the other hand, if an $\hbox{SL}(d,\Z)$-action is a restriction of an ergodic measure preserving
$\hbox{SL}(d,\R)$-action, then by a theorem of S.~Mozes (see \cite{moz}), it is mixing of all orders. 
}
\end{remark}

We conclude by proving a result of Krengel-type \cite{k}, which can be considered as a generalization of the fact that every ergodic automorphism of the torus has countable Lebesgue spectrum.

\begin{proposition}
Let $H$ be a subgroup of $\hbox{\rm SL}(2,\mathbb{Z})$ which acts in mixing fashion on $\mathbb{T}^2$.
Then for every $f\in L^2(\mathbb{T}^2)$ and every $\varepsilon >0$, there exist
$f_0\in L^2(\mathbb{T}^2)$ and a subgroup $H_0$ of finite index in $H$ such that
\begin{equation}\label{eq_kre}
\int_{\mathbb{T}^2} f_0(h\xi)\bar f_0(\xi)d\xi =0 
\end{equation}
for every $h\in H_0$, $h\ne e$.
\end{proposition}

\begin{proof}
First, we note that since $H$ is mixing, for every $(x,y)\in(\mathbb{Z}^2)^2-\{(0,0)\}$
there is at most one $h$ such that $\t hx=y$. Indeed, if $\t h_1 x =\t h_2 x$ for some $h_1,h_2\in H$,
then $h_1h_2^{-1}\in H$ is a unipotent element.

One can choose
$$
f_0=\sum_{i=1}^m a_i \chi_{_{x_i}}
$$
with some $a_i\in \mathbb{C}$ and $x_i\in\mathbb{Z}^2$ such that $\|f-f_0\|_2<\varepsilon$. We have
$$
\int_{\mathbb{T}^2} f_0(h\xi)\bar f_0(\xi)d\xi =\sum_{i,j=1}^m a_i\bar a_j \int_{\mathbb{T}^2}\chi_{_{\t hx_i -x_j}}(\xi)d\xi.
$$
Let $h_{ij}\in H-\{e\}$ be the unique element such that $\t h_{ij} x_i=x_j$ (if such an element exists).
Since $\hbox{\rm SL}(2,\mathbb{Z})$ is finitely approximable, the subgroup $H$ is finitely
approximable too. There exists a subgroup $H_0$ of finite index in $H$ such that $h_{ij}\notin H_0$
for every $i,j=1,\ldots,m$. Then (\ref{eq_kre}) holds for every $h\in H_0$, $h\ne e$.
This proves the proposition.
\end{proof}

\section*{ Acknowledgments }

We would like to thank D.~Berend, N.~Hindman, and R.~Pavlov for useful comments on the preliminary version of the paper.
We also would like to thank K.~Schmidt for bringing  Bhattacharya's paper \cite{bha} to our attention.

\end{document}